\documentclass[3p]{elsarticle}
\usepackage{amsmath,bm,accents}
\usepackage{amssymb}
\usepackage{mathrsfs}
\newtheorem{thm}{Theorem}
\newcommand{\jump}[1]{\ensuremath{[\![#1]\!]}}
\usepackage{lineno,hyperref}
\modulolinenumbers[5]
\usepackage{caption}
\journal{Journal of Computer Methods in Applied Mechanics and Engineering}

\makeatletter
\def\ps@pprintTitle{%
 \let\@oddhead\@empty
 \let\@evenhead\@empty
 \def\@oddfoot{}%
 \let\@evenfoot\@oddfoot}
\makeatother

\begin{document}

\begin{frontmatter}

\title{\small {\bf An unconditionally stable algorithm for generalised thermoelasticity based on operator-splitting and time-discontinuous Galerkin finite element methods}}


\author[a,b]{Mebratu F. Wakeni\corref{mycorrespondingauthor}}
\cortext[mycorrespondingauthor]{Corresponding author}
\ead{wknmeb001@myuct.ac.za}

\author[a,b]{ B.D. Reddy}

\author[a]{ A.T. McBride}

\address[a]{Centre for Research in Computational and Applied Mechanics (CERECAM),
University of Cape Town, 7701 Rondebosch, South Africa}
\address[b]{Department of Mathematics and Applied Mathematics,
University of Cape Town, 7701 Rondebosch, South Africa}

\begin{abstract}
An efficient time-stepping algorithm is proposed based on operator-splitting and the space-time discontinuous Galerkin finite element method for problems in the non-classical theory of thermoelasticity. The non-classical theory incorporates three models; the classical theory based on Fourier's law of heat conduction resulting in a hyperbolic-parabolic coupled system, a non-classical theory of a fully hyperbolic extension, and a combination of the two. The general problem is split into two contractive sub-problems, namely the mechanical phase and the thermal phase. Each sub-problem is  discretised using space-time discontinuous Galerkin finite element method resulting each to be stable which then leads to unconditional stability of the global product algorithm. A number of numerical examples are presented to demonstrate the performance and capability of the method.
\end{abstract}


\end{frontmatter}
\section{Introduction}
In some solids thermal energy can be transmitted by the mechanism of \emph{wave-like} propagation of heat, 
unlike the usual mechanism of conduction by diffusion. This phenomenon of heat conduction as waves, known as 
\emph{second sound}, has been observed experimentally (see, for example, \cite{Dreyer1993,Caviglia1992} for an extensive survey of experimental works involving propagation of heat as a thermal wave).  

The classical theory of heat conduction based on \emph{Fourier's law}, fails to model the second sound phenomenon. Moreover, the classical theory permits \emph{infinite speed} of propagation of parts of a localized initial heat pulse, which is paradoxical from a physical point of view. As a result, efforts have been made in an attempt to find a consistent model of heat propagation that is capable of capturing the second sound phenomenon with finite speed (see, for example, \cite{Hetnarski1999, Straughan2011} for a review of models of heat conduction as waves).

One of the alternative theory for formulating the propagation of heat in a general way that aims at 
capturing the second sound phenomenon was proposed by Green and Naghdi \cite{Green1991,Green1992,Green1993,Green1995}.
The theory of Green and Naghdi is based on three types of constitutive relations for the heat flux: Type I is equivalent to the classical theory based on Fourier's law. Type II permits the propagation of a localized heat signal as thermal wave without dissipation (see \cite{Bargmann2013} for a remark on the appropriateness of this classification). Type III is the most general theory, which includes both type I and II as special cases, in which second sound phenomenon is supported while dissipation is incorporated in the process. 

The thermomechanical coupling of non-classical heat conduction with classical elasticity is the subject of non-classical thermoelasticity. Extensive overviews of the non-classical thermoelasticity of Green and Naghdi can be found in \cite{Chandra1996, Chandra1998, Hetnarski1999}. Theoretical results concerning the non-classical theory have been addressed in several research works. In \cite{Bargmann2008c} exact solutions are obtained for thermal wave propagation in one dimension. Results on the existence and uniqueness of solutions of non-classical problem of thermoelasticity can be found, for example, in \cite{Quintanilla2002} and the references therein.  

Designing a robust and efficient numerical solution strategy for strongly coupled problems of hyperbolic-type is challenging. This is particularly the case for the non-classical 
theory of thermoelasticity where hyperbolic (or nearly hyperbolic) heat conduction equation is coupled with 
the classical 
hyperbolic elasticity problem. A standard approach for solving such time-dependent problems  is the Method 
of Lines (MoL) in which the governing partial differential equation is first discretised in space using the 
finite element method (FEM) leading to a system of ordinary differential equations, which can then be solved using the finite 
difference method. Despite its popularity, MoL struggles to accurately solve problems involving propagation 
of sharp gradients or discontinuities \cite{Hulbert1990,Hughes1988}.  

Recently, a great deal of attention has been invested in designing a spatial Discontinuous Galerkin (DG) approach for convection-dominated problems; see for example \cite{Cockburn2001}. However, these methods, 
like MoL, are based on decoupling space and time in the sense that space and time are treated differently. 
Hulbert and Hughes \cite{Hulbert1990,Hughes1988} introduced a powerful scheme based on a space-time DG finite element methodology for linear elastodynamics problems. In their approach, space and time are treated 
simultaneously and the unknown fields are allowed to be discontinuous in time while continuous in space. Recently, the space-time DG method has been used in \cite{Khalmonova2008} for classical thermoelasticity, using a monolithic approach where all the unknown fields are solved for simultaneously. 

Recently, in \cite{Bargmann2008b} a numerical solution approach based on MoL was proposed for non-classical thermoelasticity in which time integration was done in two ways: continuous Galerkin FEM for type II and III, while mixed-discountinuous Galerkin FEM for the classical problem based on the Fourier's law of heat conduction. In their approach a streamline-upwind numerical stabilization was added to localize numerical oscillations due to the propagation of sharp thermal wave. 

In the current work, we extend the existence and uniqueness results otained in \cite{Quintanilla2002} for type II theory to the more general problem of type III thermoelasticity. We also present a novel numerical algorithm for the non-classical thermoelasticity based on an operator-splitting technique motivated by Armero and Simo \cite{Armero1992} for classical thermoelasticity, coupled with a space-time DG methodology that extends the work of Hulbert and Hughes \cite{Hulbert1990} which was formulated for linear elastodynamics. The major contributions of this work are twofold: (i) the adaptation of the operator-splitting strategy for classical thermoelasticity first proposed by \cite{Armero1992} to the non-standard theory, in which the operator defining non-standard thermoelasticity is split in a way that the resulting sub-operators retain the same contractive behaviour as the global operator; and (ii) the time-DG formulation in which continuity of the unknown fields is enforced weakly by using an $L^2$-inner product in contrast to the energy-norm used in \cite{Hulbert1990}. 

The rest of this paper is organized as follows. In Section~\ref{sec:GNT}, the governing equations of the non-classical theory are summarized in a general framework of type III thermoelasticity. Well-posedness and physically meaningful boundary and initial conditions are also discussed in this section. An operator-splitting strategy for the problem of type III thermoelasticity is proposed and the resulting sub-operators are analysed in Section~\ref{sec:alg-OS}. In Section~\ref{sec:T-DG}, time-DG formulation is proposed for the sub-problems and stability of the individual algorithms and the global one is analysed in detail. A number of numerical examples both in 1--D and 2--D are presented in Section~\ref{sec:NumResults} to demonstrate the excellent performance and capability of the proposed numerical scheme. Finally, concluding remarks and some open problems are discussed in Section~\ref{sec:conclusion}. 

\section{Model problem: Non-classical thermoelasticity (type III)}\label{sec:GNT}

This section summarizes the equations governing the non-classical theory of thermoelasticty of type III as proposed by Green and Naghdi. Well-posedness of the problem is analysed. Results obtained here will serve as the basis of design and analysis of the numerical algorithm that will be presented in later sections.

\subsection*{Governing equations} 
Let $\Omega \subset\mathbb{R}^{d}$, with $1\leq d \leq 3$ be the reference placement of a continuum body 
$\mathcal{B}$ with smooth boundary $\Gamma$. Following Green and Naghdi's theory of thermoelasticity of type 
III, the system of partial differential equations governing the thermomechanical interaction in the solid 
$\mathcal{B}$ are
\begin{equation}
\left.
\begin{aligned}
\dot{\bm{u}} &= \bm{v}\label{coupled-hyperbolic1}\\
\rho \dot{\bm{v}} & = \mathrm{div}[\mathbb{C}\pmb{\bm{\varepsilon}}(\bm{u}) - \mathbf{m}\vartheta ] + \rho\mathbf{b}\\
\dot{\alpha} &= \Theta\\
\rho c \dot{\vartheta} &= \mathrm{div}[\mathbf{k}_{2}\nabla \alpha +\mathbf{k}_{3}\nabla \Theta ] - \Theta_{0}\mathbf{m}:\pmb{\bm{\varepsilon}}(\dot{\bm{u}}) + \rho r
\end{aligned}
\right\rbrace
 \text{ in } \Omega \times \mathbb{I},
\end{equation} 
where $\mathbb{I}=[0,~T]$ is the time interval of interest of length $T>0$. Superimposed dots denote time derivatives. The displacement and the velocity vector fields are denoted by $\bm{u}$ and $\bm{v}$ respectively. The scalar field $\vartheta$ denotes the relative temperature with respect to a uniform reference value $\Theta_0>0$ such that the absolute temperature $\Theta$ is given by $\Theta = \vartheta + \Theta_{0}$. The quantities $\mathbf{b}$ and $r$ are the prescribed body force and heat source.

Green and Naghdi's theory of non-classical thermoelasticity is based on the inclusion of a state variable, known as the \emph{thermal displacement} $\alpha$, that is defined in terms of an \emph{empirical temperature} $\hat{T}$ (which is assumed to coincide with the absolute temperature $\Theta$) through equation \eqref{coupled-hyperbolic1}$_{3}$, see, for example, \cite{Bargmann2008,Bargmann2008b,Bargmann2014} and the references therein.
 
The symbol $\pmb{\bm{\varepsilon}}(\bm{u})=\mathrm{sym}(\nabla \bm{u})$ denotes the small strain tensor associated with a displacement $\bm{u}$. It is assumed that the elasticity tensor $\mathbb{C}$ has the following properties:
\begin{align}
&\mathbb{C}_{ijkl} = \mathbb{C}_{jikl} = \mathbb{C}_{ijlk}, \label{eq:minor-sym}\\
&\mathbb{C}_{ijkl} = \mathbb{C}_{klij},\label{eq:major-sym}\\
&\mathbb{C}_{ijkl}\bm{\epsilon}_{ij}\bm{\epsilon}_{kl}>0 \quad \text{for any non-zero symmetric tensor }\bm{\epsilon}.\label{eq:positive-def}
\end{align}
Equations \eqref{eq:minor-sym} and \eqref{eq:major-sym} are minor and major symmetries of $\mathbb{C}$, while equation \eqref{eq:positive-def} is the positive definiteness of $\mathbb{C}$. The coupling second-order tensor $\mathbf{m}$ is of the form
\begin{equation*}\label{eq:coupling-tensor}
\mathbf{m}=3\omega\kappa\mathbf{1},
\end{equation*}
where $\omega$, $\kappa=\lambda+2/3\mu$, and $\mathbf{1}$ denote, respectively, the thermal expansion coefficient, the bulk modulus and the identity second-order tensor, and $\mu$ and $\lambda$ are the Lam\'e constants.  It is assumed that the tensor $\mathbf{k}_{2}$ is symmetric and positive-definite, and that $\mathbf{k}_{3}$ are symmetric and positive-semidefinite. The scalars $\rho>0$ and $c>0$ denote the material density and heat capacity.     
\newline
\textbf{Remarks:}
\begin{itemize}
\item[1.] The non-classical theory of thermoelasticity of type III \eqref{coupled-hyperbolic1} is the most general one in that it contains both type I and II as a special cases.  If\, $\mathbf{k}_{2}\nabla\alpha$ is omitted form \eqref{coupled-hyperbolic1}, then one obtains type I (or classical thermoelastic model) where a parabolic heat conduction equation is coupled with the hyperbolic mechanical equation. On the other hand, if\, $\mathbf{k}_{3}$ is set to zero, one obtains the type II thermoelastcity where, now, part of the system \eqref{coupled-hyperbolic1}, that is responsible for heat conduction, is hyperbolic (non-classical heat conduction).
\item[2.] Under the assumption of mechanical and thermal isotropy the elasticity tensor $\mathbb{C}$, and the tensors $\mathbf{k}_{2}$ and $\mathbf{k}_{3}$ become
\begin{equation*}
\label{eq:isotropy}
\mathbb{C} = \lambda\mathbf{1}\otimes\mathbf{1}+2\mu\mathbf{I},\quad \mathbf{k}_2 = k_2\mathbf{1},\quad\text{and}\quad  \mathbf{k}_3 = k_3\mathbf{1},
\end{equation*} 
where $\mathbf{I}$ denotes the forth-order identity tensor, and $k_2>0$, and $k_3\geq 0$ are constants. 
\item[3.] The free energy $\psi$, and hence the stress tensor $\bm{\sigma}$ and the entropy density $\eta$ are given by
\begin{equation}\label{eq:stress-entropy}
\begin{aligned}
\quad\rho \psi &= \dfrac{1}{2}\pmb{\bm{\varepsilon}}:\mathbb{C}\pmb{\bm{\varepsilon}}-\vartheta\mathbf{m}:\pmb{\bm{\varepsilon}}-\dfrac{1}{2}\dfrac{\rho c}{\Theta_0}\vartheta^2-\vartheta S_0,\\
\bm{\sigma} &= \dfrac{\partial (\rho\psi)}{\partial \pmb{\bm{\varepsilon}}}=\mathbb{C}\pmb{\bm{\varepsilon}}(\bm{u}) - \mathbf{m}\vartheta,\quad\text{and} \quad  \rho\eta =-\dfrac{\partial (\rho\psi)}{\partial \Theta}= \dfrac{c\rho}{\Theta_0}\vartheta+\mathbf{m}:\pmb{\bm{\varepsilon}}(\bm{u})+S_0,
\end{aligned}
\end{equation}
where $S_0$ is the absolute entropy density. 
\item[4.] The heat flux vector $\bm{q}$ within the non-classical theory of thermoelasticity of type III is defined as
\begin{equation*}
\bm{q} = -[\mathbf{k}_{2}\nabla \alpha +\mathbf{k}_{3}\nabla \Theta].
\end{equation*} 
\end{itemize}
Using the entropy constitutive relation \eqref{eq:stress-entropy}$_{3}$ the coupled system \eqref{coupled-hyperbolic1} can be written in terms of $\eta$ as
 \begin{equation}
\left.
\begin{aligned}
\dot{\bm{u}} &= \bm{v}\\
\rho \dot{\bm{v}} & = \mathrm{div}[\mathbb{C } \pmb{\bm{\varepsilon}} - \mathbf{m}\vartheta ] + \rho\mathbf{b}\\
\dot{\alpha} &= \Theta\\
\rho \Theta_0\dot{\eta} &= \mathrm{div}[\mathbf{k}_{2}\nabla \alpha +\mathbf{k}_{3}\nabla \Theta ] + \rho r
\end{aligned}
\right\rbrace
 \text{ in } \Omega \times [0,~T].
\end{equation} 
It is this form of the dynamical system which is crucial in designing the computational scheme based on operator-splitting in latter sections.

\subsection{Initial and boundary conditions}
Let $\{\Gamma_{\bm{u}},\Gamma_{\bm{t}}\}$ and $\{\Gamma_{\vartheta},\Gamma_{q}\}$ be two partitions of $\Gamma$, each contains mutually disjoint subsets; that is,
\begin{equation*}
 \Gamma = \overline{\Gamma_{\bm{u}} \cup \Gamma_{\bm{t}}}= \overline{\Gamma_{\vartheta} \cup \Gamma_{q}}, \text{ with } \Gamma_{\bm{u}} \cap \Gamma_{\bm{t}} =\Gamma_{\vartheta} \cap \Gamma_{q} = \emptyset.
\end{equation*}
Let $\bar{\bm{u}}:\Gamma_{\bm{u}}\times \mathbb{I}\to \mathbb{R}^{d},~$ $\bar{\bm{t}}:\Gamma_{\bm{t}}\times\mathbb{I}\to\mathbb{R}^{d},~$ $\bar{\vartheta}:\Gamma_{\vartheta}\times\mathbb{I}\to\mathbb{R},~$ and $\bar{q}:\Gamma_{q}\times\mathbb{I}\to\mathbb{R}~$ be prescribed displacement, traction, thermal displacement and flux fields. Thus the boundary conditions are given by
 \begin{equation}\label{boundary equation}
 \begin{aligned}
 \bm{u} &= \bar{\bm{u}} \quad \text{ on }\quad \Gamma_{\bm{u}}\times \mathbb{I},\\ 
 \bm{\sigma}\mathbf{n} &= \bar{\bm{t}} \quad \text{ on }\quad \Gamma_{\bm{t}}\times\mathbb{I},
 \end{aligned}
 \begin{aligned}
  \quad\qquad\vartheta &= \bar{\vartheta} \quad \text{ on } \quad \Gamma_{\vartheta}\times\mathbb{I},\\
  \qquad
\bm{q}\cdot \bm{n} &= \bar{q} \quad \text{ on } \quad \Gamma_{q}\times\mathbb{I},
\end{aligned}
 \end{equation}
 where $\bm{n}$ denotes the outward unit normal field to $\Gamma$. It is easy to observe the analogy between the two set of equations: the mechanical part \eqref{coupled-hyperbolic1}$_{1,2}$ and the thermal part \eqref{coupled-hyperbolic1}$_{3,4}$. In such analogy, we clearly see that the displacement $\bm{u}$ goes with the thermal displacement $\alpha$ (in fact, it is this analogy that motivated the name \emph{thermal displacement} \cite{Green1991}), and the velocity $\bm{v}$ goes with the absolute temperature $\Theta$, and hence with $\vartheta$. As a consequence, however, one would expect a thermal Dirichlet boundary condition is given in terms of $\alpha$ as it is customarily the case in the mechanical part that $\bm{u}$ is prescribed as Dirchlet boundary condition. The thermal Dirichlet boundary condition, in this case, is given via the relative temperature $\vartheta$ (and hence absolute temperature $\Theta$). The reason for this is that, usually, boundary conditions are prescribed in terms of physical quantities, which can be measured, which, in the thermal case, is the relative temperature, $\vartheta$ (or $\Theta$). 
 
 Furthermore, the initial conditions read 
 \begin{equation} \label{initial equation}
 \begin{aligned}
 \bm{u}(\mathbf{x},0)&=\bm{u}^{0}(\mathbf{x}),\\
  \alpha(\mathbf{x}, 0) &= \alpha^{0}(\mathbf{x}),
  \end{aligned}
  \begin{aligned}
 \quad \bm{v}(\mathbf{x}, 0) &= \bm{v}^{0}(\mathbf{x}),\\
 \quad \vartheta (\mathbf{x}, 0) &= \vartheta^{0}(\mathbf{x}),
 \end{aligned}
 \end{equation}
 where $\bm{u}^{0}$, $\bm{v}^{0}$, $\alpha^{0}$, and $\vartheta^{0}$ are prescribed initial displacement, velocity, thermal displacement and absolute temperature respectively. In prescribing an initial thermal state, a thermal configuration is assumed so that the initial thermal displacement $\alpha$ is homogeneous, that is $\alpha^0=0$, while, the physically observable quantity, the relative temperature can be initiated at a non-zero value.
  \subsection{Well-posedness: Dissipation and conservation}\label{sec:well-posed}
 Let $L_c$, $T_c$, $M_c$, and $K_c$ be characteristic scalar quantities with the dimensions of length, time, mass, and temperature, respectively. Define the dimensionless variables as
\begin{xalignat*}{4}
\bar{\bm{u}}&=\bigg[\dfrac{1}{L_c}\bigg]\mathbf{u}, & \bar{\bm{v}} &=\bigg[\dfrac{T_c}{L_c}\bigg]\bm{v}, & \bar{\bm{x}}&= \bigg[\dfrac{1}{L_c}\bigg]\bm{x}, & \bar{t} &= \bigg[\dfrac{1}{T_c}\bigg]t,\\
\bar{\Theta} &= \bigg[\dfrac{1}{K_c}\bigg]\Theta, & \bar{\alpha} &= \bigg[\dfrac{1}{T_cK_c}\bigg]\alpha, & \bar{\rho} &= \bigg[\dfrac{L_c^{3}}{M_c}\bigg]\rho, & \bar{\Theta}_0 &= \bigg[\dfrac{1}{K_c}\bigg]\Theta_0.
\end{xalignat*}
After introducing the dimensionless variables, the non-dimensional form of \eqref{coupled-hyperbolic1} become
 \begin{equation}
\left.
\begin{aligned}
\dot{\bar{\bm{u}}} &= \bar{\bm{v}},\label{nonD}\\
\rho \dot{\bar{\bm{v}}} & = \mathrm{div}[\bar{\mathbb{C}}\pmb{\bm{\varepsilon}}(\bar{\bm{u}}) - \bar{\mathbf{m}}\bar{\vartheta} ] + \bar{\rho}\bar{\mathbf{b}},\\
\dot{\bar{\alpha}} &= \bar{\Theta},\\
\bar{\rho} \bar{c} \dot{\bar{\Theta}} &= \mathrm{div}[\bar{\mathbf{k}}_{2}\nabla \bar{\alpha} +\bar{\mathbf{k}}_{3}\nabla \bar{\Theta}] - \bar{\Theta}_{0}\bar{\mathbf{m}}:\pmb{\bm{\varepsilon}}(\dot{\bar{\bm{u}}}) + \bar{\rho} \bar{r},
\end{aligned}
\right.
\end{equation}
where the spatial and time derivatives are with respect to the dimensionless space and time variables, and 
\begin{xalignat*}{4}
\bar{\mathbb{C}} &= \bigg[\dfrac{L_cT_c^2}{M_c}\bigg]\mathbb{C}, & \bar{\mathbf{m}}&=\bigg[\dfrac{L_cT_c^2K_c}{M_c}\bigg]\mathbf{m},& \bar{\bm{b}}&=\bigg[\dfrac{T_c^2}{L_c}\bigg]\bm{b},& \bar{c}&=\bigg[\dfrac{K_cT_c^2}{L_c^2}\bigg]c,\\
\bar{\mathbf{k}}_2&=\bigg[\dfrac{T_c^4K_c}{M_cL_c}\bigg]\mathbf{k}_2,&\bar{\mathbf{k}}_3&=\bigg[\dfrac{T_c^3K_c}{M_cL_c}\bigg]\mathbf{k}_3,&\bar{r}&=\bigg[\dfrac{T_c^3}{L_c^2}\bigg]r.
\end{xalignat*}
If we drop the bars in the notations of equation \eqref{nonD}, similar expressions as in equation \eqref{coupled-hyperbolic1} is obtained. For the reminder of this section, whenever the system \eqref{coupled-hyperbolic1} is mentioned, unless stated otherwise, it refers to its non-dimensional form, and the initial and boundary conditions should also be understood accordingly.

 The positive-definiteness property of $\mathbb{C}$ and $\mathbf{k}_{2}$, and the positive-semidefiniteness of $\mathbf{k}_{3}$ imply that the system \eqref{coupled-hyperbolic1} together with the initial and boundary conditions \eqref{boundary equation} and \eqref{initial equation} define an evolution equation of a general form 
 \begin{equation}\label{evolution equations}
\left.
\begin{aligned}
\dot{\bm{\chi}}(t)  &= \mathbf{A}\bm{\chi}(t) + \mathbf{f}\\
\bm{\chi}(0) & = \bm{\chi}^{0}
\end{aligned}
\right\} \quad \text{ in } \mathcal{V},
\end{equation}
where $\mathbf{A}$ is a closed linear operator with dense domain $\mathcal{D}(\mathbf{A})\subset \mathcal{V}$ defined in some suitable Banach space $\mathcal{V}$. For the case of non-classical linear thermoelasticity, for the sake of simplicity, we consider homogeneous Dirichlet boundary condition with respect to both $\bm{u}$ and $\alpha$  and the space $\mathcal{V}$  
\begin{equation}
\mathcal{V} {:=} \left\lbrace(\bm{u},\bm{v}, \alpha, \Theta)^{T}\in[\mathbf{H}^{1}(\Omega)]^{d}\times [L^2(\Omega)]^{d}\times \mathbf{H}^{1}(\Omega)\times L^2(\Omega):\bm{u}=\bm{0},\alpha = 0 \text{ on }\Gamma\right\rbrace,
\end{equation}   
is a Hilbert space.

The abstract solution vector $\bm{\chi}=(\bm{u},\bm{v}, \alpha, \Theta)^{T}\in \mathcal{V}$, while the linear operator $\mathbf{A}$ and the source term $\mathbf{f}$ in \eqref{evolution equations} are defined by 
\begin{equation}\label{eq:theOperator}
\mathbf{A}\bm{\chi}{:=}
\begin{bmatrix}
\bm{v}\\
\dfrac{1}{\rho}\mathrm{div}[\mathbb{C}\pmb{\bm{\varepsilon}}(\bm{u}) - \mathbf{m}\vartheta]\\
\Theta\\
\dfrac{1}{\rho c}\mathrm{div}[\mathbf{k}_{2}\nabla \alpha+\mathbf{k}_{3}\nabla \Theta] -\dfrac{\Theta_{0}}{\rho c}\mathbf{m}:\bm{\varepsilon}(\bm{v})
\end{bmatrix},\quad \mathbf{f}{:=}\begin{bmatrix}
\bm{0}\\
\bm{b}\\
0\\
\dfrac{1}{c}r
\end{bmatrix}.
\end{equation}
We consider an inner product, $\left\langle\cdot,\cdot\right\rangle_{\mathcal{V}}$ on $\mathcal{V}$ defined by
\begin{equation}
\left\langle\bm{\chi},\bar{\bm{\chi}}\right\rangle_{\mathcal{V}} = \left\langle\pmb{\bm{\varepsilon}}(\bm{u}),\,\mathbb{C}\pmb{\bm{\varepsilon}}(\bar{\bm{u}})\right\rangle+\left\langle\rho\bm{v},\,\bar{\bm{v}}\right\rangle + \left\langle \mathbf{k}_2^{*}\nabla \alpha,\, \nabla \bar{\alpha} \right\rangle + \left\langle c^{*}\vartheta,\,\bar{\vartheta}\right\rangle.
\end{equation}
where $\left\langle\cdot,\,\cdot\right\rangle$ denotes the standard $L^2$-inner product pairing of tensor, vector, or scalar fields, that should be understood in context, and $\mathbf{k}_{2}^{*}=\mathbf{k}_{2}\rho c/\Theta_0$ and $c^*=\rho c/\Theta_0$. The norm on $\mathcal{V}$ induced by the inner product $\left\langle\cdot,\cdot\right\rangle_{\mathcal{V}}$ is denoted by $\Vert\cdot\Vert_{\mathcal{V}}$.

Note that the linear differential operator $\mathbf{A}:\mathcal{D}(\mathbf{A})\subset \mathcal{V}\rightarrow\mathcal{V}$ is closed and the space
\begin{equation}
[\mathbf{H}^{1}_{0}(\Omega)\cap\mathbf{H}^{2}_{0}(\Omega)]^d\times[\mathbf{H}^{1}_{0}(\Omega)]^d\times(H^{1}_{0}(\Omega)\cap H^{2}_{0}(\Omega))\times H^{1}_{0}(\Omega)\subset\mathcal{D}(\mathbf{A}),
\end{equation}
is dense in $\mathcal{V}$. Hence $\mathcal{D}(\mathbf{A})$ is dense in $\mathcal{V}$.

A very important inequality concerning the evolution equation \eqref{evolution equations} is \emph{dissipativity property} of the defining operator $\mathbf{A}$. An operator $\mathbf{A}$ on closed subspace $\mathcal{D}(\mathbf{A})$ of a Hilbert space $\mathcal{V}$ endowed with an inner product $\left\langle\cdot,\cdot\right\rangle_{\mathcal{V}}$  is said to be dissipative if it satisfies the inequality $\left\langle\mathbf{A}\bm{\chi},\,\bm{\chi}\right\rangle_{\mathcal{V}}\leq 0$ for each $\bm{\chi}\in\mathcal{D}(\mathbf{A})$ \cite{Armero1992}. 
If the operator $\mathbf{A}$ is dissipative, the norm of the solution of the corresponding evolution equation is monotonically decreasing in time, which is referred to as \emph{contractivity} of the solution. That is, for a solution $\bm{\chi}$ of the evolution problem \eqref{evolution equations}, assuming dissipativity of $\mathbf{A}$ and $\mathbf{f}=\mathbf{0}$, we have 
\begin{equation}
\dfrac{d}{d t}\Vert\bm{\chi}\Vert_{\mathcal{V}}=\dfrac{d}{d t}\left\langle\bm{\chi},\bm{\chi}\right\rangle_{\mathcal{V}}=2\left\langle\dot{\bm{\chi}},\bm{\chi}\right\rangle_{\mathcal{V}}=2\left\langle\mathbf{A}\bm{\chi},\,\bm{\chi}\right\rangle\leq 0.
\end{equation}

Now, we shall show that the operator $\mathbf{A}$ that defines the problem \eqref{evolution equations} is dissipative. Let $\bm{\chi}=(\bm{u}, \bm{v}, \alpha, \vartheta)^{T}$ be in the domain of $\mathbf{A}$, $\mathcal{D}(\mathbf{A})$, satisfying the homogeneous boundary condition. Then 
\begin{align}\label{dissipative}
\nonumber \langle \bm{\chi},\mathbf{A}\bm{\chi}\rangle_{\mathcal{V}} = & \;\langle \pmb{\bm{\varepsilon}}(\bm{u}),\mathbb{C}\pmb{\bm{\varepsilon}}(\mathbf{v})\rangle + \langle \rho\mathbf{v},\dfrac{1}{\rho}\mathrm{div}[\mathbb{C}\pmb{\bm{\varepsilon}}(\bm{u})-\mathbf{m}\vartheta\rangle \\
\nonumber & + \langle \mathbf{k}_{2}^{*}\nabla \alpha, \nabla \Theta \rangle +\langle c^{*}\vartheta,\frac{1}{\rho c}\mathrm{div}[\mathbf{k}_2\nabla \alpha + \mathbf{k}_3\nabla \Theta] - \frac{\Theta_{0}}{c}\mathbf{m}:\pmb{\bm{\varepsilon}}(\bm{v})\rangle\\
\nonumber = & \;\langle \pmb{\bm{\varepsilon}}(\bm{u}),\mathbb{C}\pmb{\bm{\varepsilon}}(\bm{v})\rangle - \langle \pmb{\bm{\varepsilon}}(\bm{v}),\mathbb{C}\pmb{\bm{\varepsilon}}(\bm{u})\rangle+ \langle\pmb{\bm{\varepsilon}}(\bm{v}),\mathbf{m}\vartheta\rangle+\langle \dfrac{\rho c}{ \Theta_{0}}\mathbf{k}_2\nabla \alpha, \nabla \Theta \rangle\\
\nonumber &-\langle \dfrac{\rho c}{\Theta_{0}}\mathbf{k}_2\nabla \Theta, \nabla \alpha \rangle -\langle \dfrac{\rho c}{ \Theta_{0}}\mathbf{k}_3\nabla \Theta, \nabla \Theta \rangle - \langle \vartheta,\mathbf{m}:\pmb{\bm{\varepsilon}}(\bm{v}) \rangle \\
=& \;-\langle \dfrac{\rho c}{\Theta_{0}}\mathbf{k}_3\nabla \Theta, \nabla \Theta \rangle\leq 0. 
\end{align}
In the general context, since $\mathbf{k}_{3}$ is positive-semidefinite equation \eqref{dissipative} leads to dissipation. In the limiting case where $\mathbf{k}_3= \bm{0}$ (type II) the above argument implies conservation of energy-norm 
\begin{equation}\label{energy}
\mathscr{E}(t):=\Vert \bm{\chi}\Vert_{\mathcal{V}}^{2}=\dfrac{1}{2}\int_{\Omega}[\pmb{\bm{\varepsilon}}(\bm{u}):\mathbb{C}\pmb{\bm{\varepsilon}}(\bm{u})+\rho\bm{v}\cdot\bm{v}+\mathbf{k}_{2}^{*}\nabla\alpha\cdot\nabla\alpha + c^{*}\vartheta^2]\mathrm{d} \Omega.
\end{equation}
This is the reason why type II is also referred to as the theory of thermoelasticity\emph{ without energy dissipation}, see, for example, \cite{Green1991,Green1993,Quintanilla2002}.

Another important relation concerning the operator $\mathbf{A}$ is that it should satisfy the following: for all $\bm{\chi}^{*}\in\mathcal{V}$, there exists $\bm{\chi}$ in $\mathcal{D}(\mathbf{A})$ such that 
\begin{equation}\label{ontoness}
\bm{\chi}-\mathbf{A}\bm{\chi} = \bm{\chi}^*,
\end{equation}
in other words, the operator $(\mathbf{1}-\mathbf{A}):\mathcal{D}(\mathbf{A})\rightarrow\mathcal{V}$ is onto.

To show that $\mathbf{A}$ satisfies the relation \eqref{ontoness}, we proceed as follows: let $\bm{\chi}=(\bm{u}, \bm{v}, \alpha, \vartheta)^T$ and $\bm{\chi}^*=(\bm{u}^*, \bm{v}^*, \alpha^*, \vartheta^*)^T$ then from the definition of $\mathbf{A}$ equation \eqref{ontoness}, implies that 
\begin{equation}\label{resolvent}
\left.
\begin{aligned}
\bm{u}-\bm{v} &= \bm{u}^*,\\
\bm{v}-\dfrac{1}{\rho}\mathrm{div}[\mathbb{C}\pmb{\bm{\varepsilon}}(\bm{u})-\mathbf{m}\vartheta] &= \bm{v}^*,\\
\alpha -\vartheta &= \alpha^*,\\
\vartheta-\dfrac{1}{\rho c}\mathrm{div}[\mathbf{k}_{2}\nabla \alpha +\mathbf{k}_{3}\nabla\Theta]+\dfrac{\Theta_{0}}{c}\mathbf{m}:\pmb{\bm{\varepsilon}}(\bm{v})&= \vartheta^{*}.
\end{aligned}
\right. 
\end{equation}
Substitution of equations \eqref{resolvent}$_{1}$ and \eqref{resolvent}$_{3}$ into the remaining equations of \eqref{resolvent} leads to the (equilibrium) problem: find $\bm{\chi}=(\bm{u},\bm{v},\alpha,\vartheta)^{T}\in \mathcal{D}(\mathbf{A})$ such that $\bm{v}=\bm{u}-\bm{u}^*$, $\vartheta=\alpha-\alpha^*$ and satisfying
\begin{equation}\label{axulary}
\left.
\begin{aligned}
\rho^2\Theta_0\bm{u}-\rho\Theta_0\mathrm{div}[\mathbb{C}\pmb{\bm{\varepsilon}}(\bm{u})-\mathbf{m}\alpha] &= \accentset{\circ}{\bm{u}}, \\
\rho c\alpha-\mathrm{div}[\mathbf{k}\nabla \alpha]+\rho\Theta_{0}\mathbf{m}:\pmb{\bm{\varepsilon}}(\bm{u})&= \accentset{\circ}{\alpha},
\end{aligned}
\right.
\end{equation} 
where $\accentset{\circ}{\bm{u}}=\rho^2\Theta_0\bm{u}^{*}+\rho^2\Theta_0\bm{v}^*+\rho\Theta_0\mathrm{div}[\mathbf{m}\alpha^*]$,  $\accentset{\circ}{\alpha}=\rho c\alpha^*+\rho c\vartheta^*-\mathrm{div}[\mathbf{k}_3\nabla\alpha^*]+\rho \Theta_{0}\mathbf{m}:\pmb{\bm{\varepsilon}}(\bm{u}^*)$, and $\mathbf{k}=\mathbf{k}_2+\mathbf{k}_3$.

The weak form of equation \eqref{axulary} reads: find $\bm{\chi}=(\bm{u},\bm{v},\alpha,\vartheta)^{T}\in \mathcal{V}$ such that $\bm{v}=\bm{u}-\bm{u}^*$, $\vartheta=\alpha-\alpha^*$ and satisfying
\begin{equation}\label{weak-resolvent}
B(\bm{\chi},\bm{\xi}) = l(\bm{\xi})
\end{equation}
for all $\bm{\xi}=(\bm{w},\bm{\nu},\beta, \varpi)\in \mathcal{V}$. The bilinear form $B(\cdot,~\cdot)$ and the right hand side functional $l(\cdot)$ are given by
 
\begin{align}
\nonumber B(\bm{\chi},\bm{\xi}) &=  \langle\rho^2\Theta_0\bm{u}, \bm{w}\rangle+\langle\rho\Theta_0\mathbb{C}\pmb{\bm{\varepsilon}}(\bm{u}),\pmb{\bm{\varepsilon}}(\bm{w})\rangle-\langle\rho\Theta_0\mathbf{m}\alpha,\varepsilon(\bm{w})\rangle\\
&\quad+\langle\rho c \alpha, \beta\rangle+\langle\mathbf{k}\nabla\alpha,\nabla\beta\rangle +\langle\rho\Theta_0\mathbf{m}:\pmb{\bm{\varepsilon}}(\bm{u}),\beta\rangle, \\
l(\bm{\xi}) &= \langle\accentset{\circ}{\bm{u}},\bm{w}\rangle + \langle\accentset{\circ}{\alpha}, \beta\rangle.\label{weakresolvent-rhs}
\end{align}
Note that $\accentset{\circ}{\bm{u}}\in [\mathbf{H}^{-1}(\Omega)]^d$ and $\accentset{\circ}{\alpha}\in \mathbf{H}^{-1}(\Omega)$ and here the symbol $\langle\cdot,\cdot\rangle$ in equation \eqref{weakresolvent-rhs} represents duality pairing in their respective spaces. 

From the definition of $B(\cdot, \cdot)$, we can easily see that it is a bounded bilinear form. Since
\begin{equation}
B(\bm{\chi},\bm{\chi})=  \langle\rho^2\Theta_0\bm{u}, \bm{u}\rangle+\langle\rho\Theta_0\mathbb{C}\pmb{\bm{\varepsilon}}(\bm{u}),\pmb{\bm{\varepsilon}}(\bm{u})\rangle +\langle\rho c \alpha, \alpha\rangle+\langle\mathbf{k}\nabla\alpha,\nabla\alpha\rangle,  
\end{equation}
then $B(\cdot, \cdot)$ is $([\mathbf{H}^{1}_{0}(\Omega)]^d\times\mathbf{H}^{1}_{0}(\Omega))$-elliptic. By applying Lax-Milgram theorem we conclude that there exists $\bm{\chi}\in \mathcal{V}$ which solves the weak problem \eqref{weak-resolvent}, and hence solves equation \eqref{resolvent}. Therefore, this proves the ontoness of the resolvent operator $(\mathbf{1}-\mathbf{A})$.

In conclusion, we have seen that the operator $\mathbf{A}$ defining the non-classical linear thermoelasticity (type III) 
 
\begin{itemize}
\item[i)] is closed,
\item[ii)] has dense domain $\mathcal{D}(\mathbf{A})$ in $\mathcal{V}$,
\item[iii)] is dissipative, and 
\item[iv)] is such that $(\mathbf{1}-\mathbf{A}):\mathcal{D}(\mathbf{A})\subset\mathcal{V}\rightarrow\mathcal{V}$ is onto. 
\end{itemize}
Therefore, by the Lumer-Phillips theorem, $\mathbf{A}$ generates a strongly continuous semigroup of contractions in $\mathcal{V}$, see, for example \cite{Quintanilla2002} and the references therein. In other words, the problem \eqref{evolution equations} is well-defined and contractive. This also means that the dynamical system represented by the equation of non-classical thermoelasticity of type III is, in general, stable in the sense of Lyapunov.

\section{Algorithms based on operator-splitting strategy}\label{sec:alg-OS}
Consider an abstract evolutionary problem of the form \eqref{evolution equations}. Assume that $\mathbf{A}$ can be expressed additively as
\begin{equation}
\mathbf{A} = \mathbf{A}_1 + \mathbf{A}_2,
\end{equation}
such that the operators $\mathbf{A}_i$, $i=1,2$ define sub-problems 
\begin{equation}\label{sub-problem}
\dot{\bm{\chi}}_{_i} = \mathbf{A}_i\bm{\chi}\,; \quad \bm{\chi}_{_i}(0)=\bm{\chi}_{_i}^{0}, \quad i=1,2.
\end{equation}
Let $\mathbb{A}_{i}^{\Delta t}$ be consistent and stable time-stepping algorithms corresponding to the sub-problems \eqref{sub-problem} with $\Delta t$ representing the time step length that the algorithms step up the state vectors, $\bm{\chi}_i$, from a given time $t$ to $t+\Delta t$. A time-stepping algorithm, $\mathbb{A}^{\Delta t}$, for the global problem is obtained by taking products of the algorithms as
\begin{equation}\label{lie-trotter-kato}
\mathbb{A}^{\Delta t}=\mathbb{A}^{\Delta t}_{2}\circ\mathbb{A}^{\Delta t}_{1},
\end{equation}
The algorithm $\mathbb{A}^{\Delta t}$ is referred to as \emph{Lie-Trotter-Kato product formula} \cite{Chorin1978}. It is sometimes called a \emph{sequential split algorithm} or \emph{single pass algorithm}. In addition to the discretisation error in the individual algorithms $\mathbb{A}^{\Delta t}_{i}$, the application of the operator-splitting strategy introduces another source of error known as the \emph{splitting error}. The splitting error associated to the Lie-Trotter-Katto product formula \eqref{lie-trotter-kato} is of order of magnitude $\mathcal{O}(\Delta t)$ (see, for example, \cite{Holden2010}). It means that $\mathbb{A}^{\Delta t}$ is only first order accurate. Higher-order algorithms based on operator-splitting strategies include:
\begin{equation}\label{eq:2nd-order-splits}
\begin{aligned}
&\text{Marchuk-Strang split:}\\
&\text{Double-pass split:}
\end{aligned}
\begin{aligned}
\qquad\qquad\qquad& \mathbb{A}^{\Delta t} = \mathbb{A}^{\Delta t/2}_1\circ\mathbb{A}^{\Delta t}_2\circ\mathbb{A}^{\Delta t/2}_1\\
&\mathbb{A}^{\Delta t} =\dfrac{1}{2}(\mathbb{A}^{\Delta t}_{2}\circ\mathbb{A}^{\Delta t}_{1}+\mathbb{A}^{\Delta t}_{1}\circ\mathbb{A}^{\Delta t}_{2}). 
\end{aligned}
\end{equation} 
The Marchuk--Strang split is second-order accurate, while the double-pass split is only first-order. Note that Lie-Trotter-Kato product formula depends on the order of operations, for example in \eqref{lie-trotter-kato}, $\mathbb{A}^{\Delta t}_1$ is applied first followed by $\mathbb{A}^{\Delta t}_2$, while the order of operations does not matter in the other two operator-splitting algorithms given in \eqref{eq:2nd-order-splits}.
\subsection{Operator-splitting for non-classical thermoelasticity}
The non-classical theory of thermoelasticity is a coupling of two dynamical systems: the classical linear elasticity and the non-Fourier thermal conduction. A naive splitting of the system \eqref{coupled-hyperbolic1} into a mechanical problem under constant thermal states (isothermal) and a thermal problem with a fixed configuration will result in at most a conditionally stable algorithm even if the sub-algorithms for the two processes are unconditionally stable \cite{Armero1992}. Rather, care must be taken in splitting the two systems; this is usually dictated by an understanding of the underlying physics. In this respect, it makes sense if we split the operator $\mathbf{A}$ in \eqref{eq:theOperator} so that in the mechanical phase the entropy is held fixed (isentropic) while the temperature is allowed to vary, and in the thermal phase heat is allowed to be conducted while the configuration is fixed. In fact, in this split, it can be shown that each sub-process defines an evolution problem which is contractive, as is the global problem. Furthermore, consistent and stable algorithms for the sub-processes render a consistent and stable algorithm for the global problem by the way of operator-splitting strategy.

To this end, inspired by the work of Armero and Simo \cite{Armero1992}, taking the 4-tuple $\bm{\Sigma}=(\bm{u},\bm{v},\alpha, \eta)^{T}$ as the state variables we consider the splitting of the system of equation \eqref{coupled-hyperbolic1} into 
\begin{equation}\label{split-sys}
\left\{
\begin{aligned}
\dot{\bm{u}} &= \bm{v},\\
\rho \dot{\bm{v}} & = \mathrm{div}[\mathbb{C}\bm{\bm{\varepsilon}}(\bm{u}) - \mathbf{m}\vartheta ] + \rho 
\mathbf{b},\\
\dot{\alpha} &= 0,\\
\rho \Theta_{0} \dot{\eta} &= 0,
\end{aligned}
\right.
\quad \text{ and }\quad
\left\{
\begin{aligned}
\dot{\bm{u}} &= 0,\\
\rho \dot{\bm{v}} & = 0,\\
\dot{\alpha} &= \Theta,\\
\rho \Theta_{0} \dot{\eta}  &= \mathrm{div}[\mathbf{k}_2\nabla\alpha +\mathbf{k}_3\nabla\Theta] + \rho r.
\end{aligned}
\right.
\end{equation}
This corresponds to the additive splitting of the operator $\mathbf{A}=\mathbf{A}_1+\mathbf{A}_2$ in 
\eqref{evolution equations}
\begin{equation}\label{split-ops}
\mathbf{A}_1\bm{\Sigma}=\begin{bmatrix}
\bm{v}\\
\dfrac{1}{\rho}\mathrm{div}[\mathbb{C}\bm{\bm{\varepsilon}}(\bm{u}) - \mathbf{m}\vartheta ]\\
0\\
0
\end{bmatrix}, \hspace{2em}
\mathbf{A}_2\bm{\Sigma}=\begin{bmatrix}
\bm{0}\\
\bm{0}\\
\Theta\\
\dfrac{1}{\rho c}\mathrm{div}[\mathbf{k}_{2}\nabla \alpha+\mathbf{k}_{3}\nabla \Theta]
\end{bmatrix}.
\end{equation}
Using the same calculation as for the dissipativity of $\mathbf{A}$ in \eqref{dissipative}, we obtain the 
estimates such that for each $\bm{\Sigma}$,
\begin{equation}
\left.
\begin{aligned}
\langle \mathbf{A}_1\bm{\Sigma},\bm{\Sigma}\rangle_{\mathcal{V}} &=0,\\
 \langle \mathbf{A}_2\bm{\Sigma},\bm{\Sigma}\rangle_{\mathcal{V}} &= -\langle\dfrac{\rho c}{\Theta_0}\mathbf{k}_3\nabla\Theta, \nabla\Theta\rangle\leq 0.
\end{aligned}
\right.
\end{equation}
In general case, since $\mathbf{k}_3$ is positive-semidefinite, both operators $\mathbf{A}_1$ and 
$\mathbf{A}_2$ are dissipative. In particular, if $\mathbf{k}_3= 0$, the systems that accounts for thermal 
conduction, \eqref{split-sys}$_2$, is energy conserving, i.e. it represents heat conduction without energy 
loss in a rigid body. Moreover, it can be shown that both of the sub-operators satisfy additional conditions in order to generate strongly continuous semigroups of contraction just like how it was done in Section \ref{sec:well-posed}.

Now, having two sub-operators $\mathbf{A}_1$ and $\mathbf{A}_2$ generating contractive semigroups let us assume that there corresponds two algorithms $\mathbb{A}^{\Delta t}_i$, $i=1,2$ which are \emph{B-stable} (\emph{non-linearly stable}) ; that is 
\begin{equation}
\Vert \mathbb{A}^{\Delta t}_{i}\bm{\chi}^{n}-\mathbb{A}^{\Delta t}_{i}\tilde{\bm{\chi}}\Vert\leq \Vert \bm{\chi}^{n}-\tilde{\bm{\chi}}\Vert \quad \mathrm{for~all~}\bm{\chi}^{n}, ~\tilde{\bm{\chi}}\in \mathcal{V}.
\end{equation}
 In the linear case this means that each of the algorithms satisfy the estimate (see \cite{Simo1991, Armero1992} and the references therein):
\begin{equation}
\Vert\mathbb{A}^{\Delta t}_{i}\Vert_{\mathcal{V}^{*}}\leq 1,\;\qquad i=1, 2.
\end{equation}  
where $\mathcal{V}^*$ is the dual of $\mathcal{V}$. Let $\left\{\bm{\chi}^{n}\right\}_{n\in\mathbb{N}}$ and $\left\{\tilde{\bm{\chi}}^{n}\right\}_{n\in\mathbb{N}}$ be two sequences in $\mathcal{V}$ generated by the Lie-Trotter-Kato product formula $\mathbb{A}^{\Delta t}$ corresponding to two initial conditions $\bm{\chi}(0) = \bm{\chi}^{0}$ and $\tilde{\bm{\chi}}(0) = \tilde{\bm{\chi}}^{0}$ respectively. Then the product formula $\mathbb{A}^{\Delta t}$ satisfies  the stability estimate
\begin{align}
\nonumber\Vert\bm{\chi}^{n+1}-\tilde{\bm{\chi}}^{n+1}\Vert_{\mathcal{V}}
\nonumber&=\Vert\mathbb{A}^{\Delta t}\bm{\chi}^{n}-\mathbb{A}^{\Delta t}\tilde{\bm{\chi}}^{n}\Vert_{\mathbf{V}}\\
\nonumber&=\Vert\mathbb{A}^{\Delta t}_2[\mathbb{A}^{\Delta t}_1\bm{\chi}^{n}]-\mathbb{A}^{\Delta t}_2[\mathbb{A}^{\Delta t}_1\tilde{\bm{\chi}}^{n}\Vert_{\mathcal{V}}\\
\nonumber&\leq \Vert\mathbb{A}^{\Delta t}_1\bm{\chi}^{n}-\mathbb{A}^{\Delta t}_1\tilde{\bm{\chi}}^{n}\Vert_{\mathcal{V}}\qquad\qquad\text{($\Vert\mathbb{A}^{\Delta t}_2\Vert_{\mathcal{V}}\leq 1$)}\\
&=\Vert\bm{\chi}^{n}-\tilde{\bm{\chi}}^{n}\Vert_{\mathcal{V}}\;\,\quad\qquad\qquad\qquad\text{($\Vert\mathbb{A}^{\Delta t}_1\Vert_{\mathcal{V}}\leq 1$)},
\end{align}
which proves the \emph{non-linear stability} of the global algorithm corresponding to the product formula $\mathbb{A}^{\Delta t}$.

The numerical scheme that is going to be formulated in the subsequent sections is based on the operator-splitting approach and time-discontinuous Galerkin finite element method. A Lie-Trotter-Kato product formula is applied to merge algorithms for the two phases. Hence, the product formula can also be viewed in the sense of \emph{predictor--corrector} regime, where the sub-algorithm for the mechanical phase is used as predictor and that of the thermal phase as a corrector. 

\section{Time-discontinuous Galerkin finite element method (T-DG FEM)}\label{sec:T-DG}
Let $\bm{T}_{h}=\left\{\Omega^e\right\}$ be a triangulation of $\bar{\Omega}$, where $\bar{\Omega}$ denotes the closure (the union of the interior and boundary) of $\Omega$, such that
\begin{equation}
\bar{\Omega} = \bigcup_{\Omega^e\in\bm{T}_{h}}\Omega^e.
\end{equation}
Denote the space of scalar piecewise polynomials on the mesh $\bm{T}_{h}$ by $\mathscr{P}_{h}^j$:
\begin{equation}
\mathscr{P}_{h}^j = \left\{\varphi_{_h}\in C^0(\bar{\Omega}): \varphi_{_h}\vert_{_{\Omega^e}}\in P^{j}
(\Omega^e),\;\Omega^e\in\bm{T}_h\right\},
\end{equation}
where $P^j(\Omega^e)$ denotes the set of polynomials of degree at most $j$ defined on $\Omega^e$.

Consider a partition of the time domain $\mathbb{I}=[0,T]$ into the collection $\left\{I_n=[t_n, t_{n+1}]\right\}_{n=0}^{N-1}$, $N\in\mathbb{N}$ of non-overlapping subintervals. The time step length is ${\Delta t}_{n}=t_{n+1}-t_n$ for $n=0,1,\cdots,N-1$ with
\begin{equation*}
0=t_0<t_1<\cdots<t_{N}=T.
\end{equation*}
For each time sub-domain $I_n$ we consider the space-time domain of the form $Q_n=\bar{\Omega}\times I_n$ 
referred to as the $n^{\text{th}}$ \emph{space-time slab}. 

Admissible scalar functions, $\phi^{h}$, that we consider in T-DG FEM will be polynomials in time $t\in I_n$ 
with coefficients from the spatial function space $\mathscr{P}^{j}_{h}$; i.e.
\begin{equation}
\phi^{h}(\bm{x},t) = \sum_{i}\varphi^{h}_{i}(\bm{x})t^{i},\hspace{1em}\varphi^{h}_{i}\in\mathscr{P}^{j}_{h},
\end{equation} 
 where $t^{i}$ is a monomial in $t\in I_n$ of order $i\in\mathbb{N}$. Denote by $\mathcal{S}_{_{h}}(Q_n;j,l)$, the space of admissible functions on the space-time domain $Q_n$ of degree $j+l$ (that is, $j$ in space and $l$ in time):
\begin{equation}
\mathcal{S}_{_{h}}(Q_n;j,l) = \left\{\phi^{h}:\phi^{h}(\bm{x}, t)=\sum_{i=0}^{l}\varphi^{h}_{i}
(\bm{x})t^{i},\;\varphi^{h}_{i}\in\mathscr{P}^{j}_{h},\;(\bm{x},t)\in Q_n\right\}.
\end{equation} 
In fact, it is easy to observe that the space $\mathcal{S}_{_{h}}(Q_n;j,l)$ is generated by the tensor 
products of the basis elements of the spaces $\mathscr{P}_{h}^j$ and $P^{l}(I_n)$--the set of polynomials in time of degree at most $j$.
\subsection*{Remark:}
\begin{itemize}
\item[1.] The space-time mesh for $Q_n$ is composed of cells with one element thickness in the time direction; and in each cell in the slab, time and space are orthogonal to each other. i.e. each cell is of the form $\Omega^e\times I_n$. Nevertheless, the formulation being developed here can be easily modified in terms of non-orthogonal space-time elements and with slabs composed of more than one element thickness in time direction as well.
\item[2.] The approach would readily accommodate the use of an adaptive mesh refinement procedures. In such cases, there may be cells with time direction thickness less than ${\Delta t}_{n}$ embedded in slab $Q_n$. In this case, at each hanging node the solution must be constrained so that the hanging nodes will be condensed out later. 
\end{itemize}  
\subsection*{Notations} Let $\varphi$ and $\psi$ be functions defined on space time domain slab $Q_n$. Some frequently used notations are 
\begin{itemize}
\item[a)] Spatial $L^2$ inner product
\begin{equation*}
\langle \varphi, \psi\rangle {:=} 
\int_{\Omega}\varphi\psi~\mathrm{d} \Omega.
\end{equation*}
\item[b)] Space-time $L^2$ inner product 
\begin{equation*}
(\varphi, \psi)_{Q_n} {:=} \int_{I_{n}} \langle \varphi, 
\psi\rangle ~\mathrm{d}t.
\end{equation*}
\item[c)] Space-time boundary integrals on $Z_n=\Gamma_{\mathbf{t}}\times I_n$ and $F_n=\Gamma_{q}\times I_n$
\begin{equation*}
\begin{aligned}
(\varphi, \psi)_{Z_n}&=\int_{I_n}\int_{\Gamma_{\mathbf{t}}}\varphi\psi~\mathrm{d}\Gamma~\mathrm{d}t,\\
(\varphi, \psi)_{F_n}&=\int_{I_n}\int_{\Gamma_{q}}\varphi\psi~\mathrm{d}\Gamma~\mathrm{d}t.
\end{aligned}
\end{equation*}
\item[c)] Right/left limit of a discontinuous function in time at $t_n$
\begin{equation*}
\varphi(t_{n}^{\pm}){:=}\lim_{\varepsilon \to 0^{\pm}} \varphi(t_{n}+\varepsilon).
\end{equation*}
\item[d)] Temporal jump of a discontinuous function at $t_n$
\begin{equation*}
\jump{\varphi}_{n} = \varphi(t_{n}^{+})-\varphi(t_{n}^{-}).
\end{equation*}
\end{itemize}

\subsection{Mechanical problem}

In the mechanical phase the entropy is held fixed, so that the last equation of the left hand system in \eqref{split-sys}; that is the equation 
\begin{equation}
\rho \Theta_{0} \dot{\eta}=\dfrac{\mathrm{d}}{\mathrm{d} t}[\rho c \vartheta + \rho \Theta_{0}\mathbf{m}:\bm{\varepsilon}(\bm{u})+S_0] = 0\hspace{1em}\text{ on } \quad\Omega\times [t_n,t_{n+1})
\end{equation} 
is solved in closed form to obtain an \emph{intermediate temperature} $\vartheta^I$. 
This leads to an explicit formula for $\vartheta^I$ in terms of the state variables at time step $t_n$ from the left and at time value $t\in (t_n,~t_{n+1})$, that is, 
\begin{equation}
\vartheta^I(t) = \vartheta(t_n^{-})-\dfrac{\Theta_0}{c}\mathbf{m}:\bm{\varepsilon}(\bm{u}(t)-\bm{u}(t_n^-))
\hspace{2em}\text{for }\;t\in (t_n, t_{n+1}).
\end{equation}
We substitute this result into the mechanical problem \eqref{split-sys}$_{1}$ to obtain
\begin{equation}
\left.
\begin{aligned}\label{mech-prob}
\dot{\bm{u}} &= \bm{v},\\
\rho \dot{\bm{v}} & = \mathrm{div}[\mathbb{C}_{_{ad}}\bm{\bm{\varepsilon}}(\bm{u})] + \bm{f},
\end{aligned}
\right.
\end{equation}
where $\mathbb{C}_{_{ad}}=\mathbb{C}+(\Theta_0/c)\mathbf{m}\otimes\mathbf{m}$, in the terminology used in \cite{Armero1992}, is referred to as \emph{adiabatic elasticity tensor} and $\bm{f} = \rho \bm{b}-\mathbf{m}[\vartheta(t_n^-)+(\Theta_0/c)\mathbf{m}:\bm{u}(t_n^-)]$ and $\vartheta(t_n^-)$ and $\bm{u}(t_n^-)$ denote the temperature at the end of the previous space-time slab, $Q_{n-1}$. Note that the adiabatic elasticity tensor $\mathbb{C}_{_{ad}}$ remains positive-definite and symmetric as $\mathbb{C}$. 

On the current space-time slab, $Q_n$, we use same boundary conditions as given in \eqref{boundary equation} for the mechanical fields but the initial conditions, in this case, are the solution for $\bm{u}$ and $\bm{v}$ at the end of the previous slab; that is $\bm{u}(t_n^{-})$ and $\bm{v}(t_n^{-})$. 

To define the T-DG FEM formulation of the mechanical problem \eqref{mech-prob}, we first define the trial and weight function spaces for displacement $\bm{u}$ and velocity $\bm{v}$ vector fields as
\begin{equation}
\begin{aligned}
\mathcal{T}^{^u}_{_h} &= \left\{\bm{u}^{h}\in[\mathcal{S}_{_{h}}(Q_n;j,l)]^{d}:\bm{u}^{h}=\bar{\bm{u}} \text{ on }\Gamma_{\bm{u}}\times I_n\right\},\\
\mathcal{T}^{^v}_{_h} &= \left\{\bm{v}^{h}\in[\mathcal{S}_{_{h}}(Q_n;j,l)]^{d}:\bm{v}^{h}=\dot{\bar{\bm{u}}} \text{ on }\Gamma_{\bm{u}}\times I_n\right\},\\
\mathcal{W}^{^u}_{_h} &= \left\{\bm{w}^{h}\in[\mathcal{S}_{_{h}}(Q_n;j,l)]^{d}:\bm{w}^{h}=\bm{0} \text{ on }\Gamma_{\bm{u}}\times I_n\right\},\\
\mathcal{W}^{^v}_{_h} &= \left\{\bm{\varphi}^{h}\in[\mathcal{S}_{_{h}}(Q_n;j,l)]^{d}:\bm{\varphi}^{h}=\bm{0} \text{ on }\Gamma_{\bm{u}}\times I_n\right\},
\end{aligned}
\end{equation}
where $\mathcal{T}^{^u}_{_h}$ and $\mathcal{T}^{^v}_{_h}$, $\mathcal{W}^{^u}_{_h}$ and $\mathcal{W}^{^v}_{_h}$ are trial and weight function spaces for displacement and velocity vector fields respectively. The T-DG FEM is formulated as: find $\bm{U}^{h}=(\bm{u}^{h},\bm{v}^{h})^{T}\in \mathcal{T}^{^u}_{_h}\times\mathcal{T}^{^v}_{_h}$ such that for all $\bm{V}^{h}=(\bm{w}^{h},\bm{\varphi}^{h})^{T}\in \mathcal{W}^{^u}_{_h}\times \mathcal{W}^{^v}_{_h}$
\begin{equation}\label{M-DG-Form}
A_{n}^{^M}(\bm{U}^{h},\bm{V}^{h}) = b_{n}^{^M}(\bm{V}^{h}),
\end{equation} 
where
\begin{equation}
\begin{aligned}
A_{n}^{^M}(\bm{U}^{h},\bm{V}^{h}) &= (\dot{\bm{u}}^{h},\bm{w}^{h})_{Q_n}-(\bm{v}^{h},\bm{w}^{h})_{Q_n}+(\rho\dot{\bm{v}}^{h},\bm{\varphi}^{h})_{Q_n}+(\mathbb{C}_{_{ad}}\bm{\varepsilon}(\bm{u}^{h}),\bm{\varepsilon}(\varphi))_{Q_n}\\
&\qquad +\langle\bm{u}^{h}(t_n^+),\bm{w}^{h}(t_n^+)\rangle + \langle\rho\bm{v}^{h}(t_n^+),\bm{\varphi}^{h}(t_n^+)\rangle\\
b_{n}^{^M}(\bm{V}^{h}) &= (\bar{\mathbf{t}},\bm{\varphi}^{h})_{Z_n} + (\bm{f},\bm{\varphi}^{h})_{Q_n} + \langle\bm{u}^{h}(t_n^-),\bm{w}^{h}(t_n^+)\rangle + \langle\rho\bm{v}^{h}(t_n^-),\bm{\varphi}^{h}(t_n^+)\rangle.
\end{aligned}
\end{equation}
\subsubsection*{Remark:}
\begin{itemize}
\item[(1)] The main difference between the DG formulation presented here in \eqref{M-DG-Form} and that of \cite{Hulbert1990} is the inner product used to enforce the equation of motion \eqref{mech-prob} weakly. In our formulation we use the $L^2$-inner product to weakly enforce the mechanical problem while in \cite{Hulbert1990} an energy-inner product is used.
\item[(2)] The formulation \eqref{M-DG-Form} is consistent in the sense of a time-stepping algorithm. This can be seen from the Euler-Lagrange form of \eqref{M-DG-Form} given by 
\begin{align}\label{eq:M-EulLag}
\nonumber 0&= A_{n}^{^M}(\bm{U}^{h},\bm{V}^{h}) - b_{n}^{^M}(\bm{V}^{h})\\
\nonumber  &= (\dot{\bm{u}}^h-\bm{v}^h,\bm{w}^h)_{Q_n}+(\rho\dot{\bm{v}}^h-\mathrm{div}[\mathbb{C}_{_{ad}}\bm{\varepsilon}(\bm{u})^h]-\bm{f},\bm{\varphi}^h) \hspace{2em}\text{ (equation of motion)}\\
\nonumber  &\quad+\langle\jump{\bm{u}^{h}}_{n},\bm{w}^{h}(t_n^+)\rangle\hspace{16em}\text{(displacement continuity)}\\
  &\quad+ \langle\jump{\rho\bm{v}^{h}}_{n},\bm{\varphi}^{h}(t_n^+)\rangle ,\hspace{14em}\,\;\;\;\;\text{(velocity continuity)}
\end{align} 
that upon substitution of a sufficiently smooth solution pair $(\bm{u},\bm{v})^{T}$ of the strong form \eqref{mech-prob} into \eqref{eq:M-EulLag}, the weak forms of the jumps and the equation of motion vanish.  
\item[(3)] The jump terms are used to improve the stability of the scheme without degrading the accuracy. As a result, the formalism can be readily extended to the non-linear case without eliminating the jump term from the displacement-velocity relation.
\item[(4)] One of the consequences of using the $L^2$-inner product is that a Dirchlet-type boundary condition may not be necessary to define the velocity trial and weight function spaces. Instead we can use
\begin{equation}
\mathcal{T}^{^v}_{_h}=\mathcal{W}^{^v}_{_h}=[\mathcal{S}_{_{h}}(Q_n;j,l)]^{d}.
\end{equation}
\end{itemize}
\subsubsection{Stability: The mechanical algorithm}
For the sake of simplicity, we assume a homogeneous source term $\bm{f}=\bm{0}$ and boundary conditions, i.e. $\bar{\mathbf{t}}=\bm{0}$ and $\bar{\bm{u}}=\bm{0}$. We claim that the formulation \eqref{M-DG-Form} renders an unconditionally stable time-stepping algorithm. That is to say:
\begin{equation}\label{eq:energy-estimate}
\mathscr{E}_{_{M}}(\bm{U}^{h}(t_{n+1}^{-}))\leq\mathscr{E}_{_M}(\bm{U}^{h}(t_{n}^{-})) \quad \forall n=0,1,\cdots,N-1,
\end{equation}
where $\mathscr{E}_{_{M}}(\bm{U}(t))$ is the total mechanical energy of $\bm{U}=(\bm{u},\bm{v})^T$ at time $t$, given  by
\begin{equation}
\mathscr{E}_{_{M}}(\bm{U}(t))=\dfrac{1}{2}\int_{\Omega}\big[\bm{\varepsilon}(\bm{u}(t)):\mathbb{C}_{_{ad}}\bm{\varepsilon}(\bm{u}(t))+\rho\bm{v}(t)\cdot\bm{v}(t)\big]~\mathrm{d}\Omega.
\end{equation}
For the analysis we use elliptic and $L^2$ interpolation operators $\bm{\pi}:[H_{\Gamma_{u}}(\Omega)]^d\rightarrow \mathcal{W}^{^u}_{_h}(t)$ and $\hat{\bm{\pi}}:[L^{2}(\Omega)]^d\rightarrow\mathcal{W}^{^u}_{_h}(t)$ respectively defined as: for $\bm{u}\in [H_{\Gamma_{u}}(\Omega)]^d$ and $\bm{w}\in [L^{2}(\Omega)]^d$
\begin{equation}\label{proj-operators1}
\begin{aligned}
\langle\mathbf{C}_{_{ad}}\bm{\varepsilon}(\bm{\pi}\bm{u})],\bm{\varepsilon}(\bm{\varphi}^h)\rangle &= \langle\mathbf{C}_{_{ad}}\bm{\varepsilon}(\bm{u}), \bm{\varepsilon}(\bm{\varphi}^h)\rangle,\qquad\forall \bm{\varphi}^h\in\mathcal{W}^{^u}_{_h}(t),\\
\langle\hat{\bm{\pi}}\bm{w},\bm{\psi}^h\rangle&= \langle\bm{w},\bm{\psi}^h\rangle, \qquad \hspace{4em}\forall\bm{\psi}^h\in\mathcal{W}^{^u}_{_h}(t),
\end{aligned}
\end{equation} 
where $\mathcal{W}^{^u}_{_h}(t)$ referees to the space of functions in $\mathcal{W}^{^u}_{_h}$ at a fixed but arbitrary $t\in I_n$. We also use the fact that \cite{Johnson1993} 
\begin{equation}\label{H-to-L}
\mathrm{div}[\mathbb{C}_{_{ad}}\bm{\varepsilon}(\bm{\pi}\bm{u})] = \hat{\bm{\pi}}\mathrm{div}[\mathbb{C}_{_{ad}}\bm{\varepsilon}(\bm{u})].
\end{equation}
Now, given the solution $\bm{U}^{h}(t_n^-)=(\bm{u}^{h}(t_{n}^-),\bm{v}^{h}(t_{n}^-))^{T}$ of \eqref{M-DG-Form} at the end of the previous space-time slab, $Q_{n-1}$, and let $\bm{U}^{h}(t)=(\bm{u}^{h}(t),\bm{v}^{h}(t))^{T}$ be the solution of \eqref{M-DG-Form} in the current space-time slab, $Q_{n}$. Replace $\bm{V}^{h}=(\hat{\bm{\pi}}\mathrm{div}[\mathbb{C}_{_{ad}}\bm{\varepsilon}(\bm{u}^h)],\bm{0})^T$ in \eqref{M-DG-Form} to obtain
\begin{equation}\label{pL-on-WF}
\begin{aligned}
0&=(\dot{\bm{u}}^{h},\hat{\bm{\pi}}\mathrm{div}[\mathbb{C}_{_{ad}}\bm{\varepsilon}(\bm{u}^h)])_{Q_n}-(\bm{v}^{h},\hat{\bm{\pi}}\mathrm{div}[\mathbb{C}_{_{ad}}\bm{\varepsilon}(\bm{u}^h)])_{Q_n}\\
&\quad +\langle\jump{\bm{u}^{h}}_n,\hat{\bm{\pi}}\mathrm{div}[\mathbb{C}_{_{ad}}\bm{\varepsilon}(\bm{u}^h(t_n^+))]\rangle,
\end{aligned}
\end{equation}
Applying \eqref{H-to-L} in \eqref{pL-on-WF} and the definition of the projection operator $\bm{\pi}$ and using integration by parts (note the homogeneous boundary conditions) leads to
\begin{equation}\label{u-to-v}
\begin{aligned}
0&=(\bm{\varepsilon}(\dot{\bm{u}}^{h}),\mathbb{C}_{_{ad}}\bm{\varepsilon}(\bm{u}^h))_{Q_n}-(\bm{\varepsilon}(\bm{v}^{h}),\mathbb{C}_{_{ad}}\bm{\varepsilon}(\bm{u}^h))_{Q_n}\\
&\quad +\langle\jump{\bm{\varepsilon}(\bm{u}^{h})}_n,\mathbb{C}_{_{ad}}\bm{\varepsilon}(\bm{u}^h(t_n^+))\rangle.
\end{aligned}
\end{equation}
Again substituting $\bm{V}^{h}=(\bm{0},\bm{v}^h)^T$ into \eqref{M-DG-Form} yields
\begin{equation}\label{momentum}
0=(\rho\dot{\bm{v}}^{h},\bm{v}^{h})_{Q_n}+(\mathbb{C}_{_{ad}}\bm{\varepsilon}(\bm{u}^{h}),\bm{\varepsilon}(v))_{Q_n}+\langle\jump{\rho\bm{v}^h}_n,\bm{v}^h(t_{n}^+)\rangle.
\end{equation}
Adding the equations \eqref{u-to-v} and \eqref{momentum} we obtain
\begin{equation}\label{eq:full}
(\bm{\varepsilon}(\dot{\bm{u}}^{h}),\mathbb{C}_{_{ad}}\bm{\varepsilon}(\bm{u}^h))_{Q_n}+(\rho\dot{\bm{v}}^{h},\bm{v}^{h})_{Q_n}+\langle\jump{\bm{\varepsilon}(\bm{u}^{h})}_n,\mathbb{C}_{_{ad}}\bm{\varepsilon}(\bm{u}^h(t_n^+))\rangle+\langle\jump{\rho\bm{v}^h}_n,\bm{v}^h(t_{n}^+)\rangle=0.
\end{equation}
Taking the time derivative out of the space integral, \eqref{eq:full} becomes 
\begin{equation}
\begin{aligned}
&\dfrac{1}{2}(\bm{\varepsilon}(\bm{u}^{h}(t_{n+1}^{-})),\mathbb{C}_{_{ad}}\bm{\varepsilon}(\bm{u}^h(t_{n+1}^-)))_{Q_n}+\dfrac{1}{2}(\rho\bm{v}^{h}(t_{n+1}^-),\bm{v}^{h}(t_{n+1}^{-}))_{Q_n}\\
&-\dfrac{1}{2}(\bm{\varepsilon}(\bm{u}^{h}(t_{n}^{+})),\mathbb{C}_{_{ad}}\bm{\varepsilon}(\bm{u}^h(t_{n}^+)))_{Q_n}-\dfrac{1}{2}(\rho\bm{v}^{h}(t_{n}^+),\bm{v}^{h}(t_{n}^{+}))_{Q_n}\\
&+\langle\jump{\bm{\varepsilon}(\bm{u}^{h})}_n,\mathbb{C}_{_{ad}}\bm{\varepsilon}(\bm{u}^h(t_n^+))\rangle+\langle\jump{\rho\bm{v}^h}_n,\bm{v}^h(t_{n}^+)\rangle=0.
\end{aligned}
\end{equation}
After some algebraic manipulation we obtain
\begin{equation}\label{energy-estimate}
\mathscr{E}_{_{M}}(\bm{U}^{h}(t_{n+1}^{-}))+\mathscr{E}_{_{M}}(\jump{\bm{U}^h}_n)=\mathscr{E}_{_M}(\bm{U}^{h}(t_{n}^{-})).
\end{equation}
Since $\mathscr{E}_{_{M}}(\jump{\bm{U}^h}_n)$ is non-negative the energy equation \eqref{energy-estimate} leads to the estimate \eqref{eq:energy-estimate} that renders the scheme for the mechanical phase \eqref{M-DG-Form} unconditionally stable. In fact, the total numerical dissipation added is precisely equal to
\begin{equation*}
\sum_{n=0}^{N-1}\mathscr{E}_{_{M}}(\jump{\bm{U}^h}_n).
\end{equation*}
\subsubsection*{Remark:}
The use of projection operators in \eqref{proj-operators1} and \eqref{pL-on-WF} reveals an important point, that is, the current T-DG formulation can be converted into the one given in \cite{Hulbert1990}.   
\subsection{Thermal problem}
The solution $\bm{U}^{h}(t_{n+1}^{-})$ of the mechanical phase is known at the end of the current space-time slab. The objective, in the present phase, is to solve for the thermal states $\bm{\Pi}(t)=(\alpha(t),\vartheta(t))^{T}$, $t\in I_n$.  Hence, the global solution at the end of the current slab will be $(\bm{U}^{h}(t_{n+1}^{-}),\bm{\Pi}^h(t_{n+1}^-))$. 

The operator-splitting is performed based on the state vector $\bm{\Sigma}=(\bm{u},\bm{v},\alpha, \eta)^{T}$. As a result, we enforce the problem in the thermal phase using the conservation form

\begin{equation}\label{T-Strong}
\begin{aligned}
\dot{\alpha} &= \Theta,\\
\rho \Theta_{0} \dot{\eta}  &= \mathrm{div}[\mathbf{k}_2\nabla\alpha + \mathbf{k}_3\nabla\Theta] + \rho r.
\end{aligned}
\end{equation}
Here, recall that the entropy density $\eta$ is obtained from the relation $\rho\eta = \dfrac{c\rho}{\Theta_0}\vartheta+\mathbf{m}:\pmb{\bm{\varepsilon}}(\bm{u})+S_0$. 

In the thermal phase, the displacement $\bm{u}$ and the velocity $\bm{v}$ are fixed at the corresponding values at the end of the current slab in the mechanical phase; that is,
\begin{equation}
\bm{u}_{T}(t)=\bm{u}_{M}(t_{n+1}^{-}), \quad \bm{v}_{T}(t)=\bm{v}_{M}(t_{n+1}^{-}), \qquad t\in I_n, 
\end{equation}
where the subscripts $T$ and $M$ represents the values of the fields in the thermal and mechanical phases, respectively. 
Thus, the time derivative of the terms $\mathbf{m}:\pmb{\bm{\varepsilon}}(\bm{u})+S_0$ vanishes, and consequently the left hand side of equation \eqref{T-Strong}$_{2}$ becomes
\begin{equation}\label{entropy-rate}
\rho\Theta_0\dot{\eta} = \rho c \dot{\vartheta}.
\end{equation}  

In addition, the jump in the entropy, in this case, reads as
\begin{equation}\label{entropy-jump}
\begin{aligned}
\jump{\rho\Theta_0\eta}_{n} &= \rho\Theta_0\eta(t_n^+)-\rho\Theta_0\eta(t_n^-)\\
                        &= \big[\rho c \vartheta(t_{n}^+)+\rho\Theta_0\mathbf{m}:\bm{\varepsilon}(\bm{u}_{M}(t_{n+1}^{-}))+S_0\big]-\big[\rho c \vartheta(t_{n}^{-})+\rho\Theta_0\mathbf{m}:\bm{\varepsilon}(\bm{u}(t_{n}^{-}))+S_0\big]\\
                       &= \jump{\rho c \vartheta}_{n}+\rho\Theta_0\mathbf{m}:[\bm{\varepsilon}(\bm{u}_{M}(t_{n+1}^-))-\bm{\varepsilon}(\bm{u}(t_n^-))].  
\end{aligned}
\end{equation} 
It should not cause any confusion if we drop the superscript $M$ in the equation \eqref{entropy-jump} so that the jump term can be written as
\begin{equation}
\jump{\rho\Theta_0\eta}_{n} =  \jump{\rho c \vartheta}_{n}+\rho\Theta_0\mathbf{m}:[\bm{\varepsilon}(\bm{u}(t_{n+1}^-))-\bm{\varepsilon}(\bm{u}(t_n^-))].
\end{equation}

To define the T-DG formulation for the thermal phase we first define the thermal displacement and temperature trial and weight function spaces $\mathcal{T}^{^\alpha}_{h}$, $\mathcal{T}^{^\vartheta}_{h}$ and $\mathcal{W}^{^\alpha}_{h}$, $\mathcal{W}^{^\vartheta}_{h}$  respectively based on $\mathcal{S}_{_{h}}(Q_n;j,l)$ and the boundary condition requirements.
\begin{equation}
\begin{aligned}
\mathcal{T}^{^\alpha}_{_h} &= \left\{\alpha^{h}\in\mathcal{S}_{_{h}}(Q_n;j,l):\dot{\alpha}^{h}=\bar{\Theta} \text{ on }\Gamma_{\alpha}\times I_n\right\},\\
\mathcal{T}^{^\vartheta}_{_h} &= \left\{\vartheta^{h}\in\mathcal{S}_{_{h}}(Q_n;j,l):\Theta^{h}=\bar{\Theta} \text{ on }\Gamma_{\alpha}\times I_n\right\},\\
\mathcal{W}^{^\alpha}_{_h} &= \left\{\beta^{h}\in\mathcal{S}_{_{h}}(Q_n;j,l):\beta^{h}=0 \text{ on }\Gamma_{\alpha}\times I_n\right\},\\
\mathcal{W}^{^\vartheta}_{_h} &= \left\{\sigma^{h}\in\mathcal{S}_{_{h}}(Q_n;j,l):\sigma^{h}=0 \text{ on }\Gamma_{\alpha}\times I_n\right\}.\\
\end{aligned}
\end{equation}
Formally, the T-DG FEM formulation of the thermal phase on the domain $Q_n$ is defined as: find $\bm{\Pi}^h=(\alpha,\vartheta)^{T}\in\mathcal{T}^{^\alpha}_{_h}\times \mathcal{T}^{^\vartheta}_{_h}$ such that for each $\bm{\Lambda}^{h}=(\beta^{h},\sigma^{h})^{T}\in\mathcal{W}^{^\alpha}_{_h}\times\mathcal{W}^{^\vartheta}_{_h}$
\begin{equation}\label{DG-T-Form}
A_{_n}^{^T}(\bm{\Pi}^h,\bm{\Lambda}^{h})=b_{_n}^{^T}(\bm{\Lambda}^{h}),
\end{equation}
where
\begin{equation}
\begin{aligned}
A_{_n}^{^T}(\bm{\Pi}^h,\bm{\Lambda}^{h}) &= (\dot{\alpha}^h, \beta^h)_{Q_n}-(\Theta^h, \beta^h)_{Q_n}+(\rho c \dot{\vartheta}^h, \sigma^h)_{Q_n} + ([\mathbf{k}_2\nabla\alpha^h+\mathbf{k}_3\nabla \theta^h], \nabla \sigma^h)_{Q_n},\\
&\quad + \langle\alpha^h(t_n^+),\beta^h(t_n^+)\rangle + \langle\rho c \vartheta^h(t_n^+), \sigma^{h}(t_n^+)\rangle\\
b_{_n}^{^T}(\bm{\Lambda}^{h})&=\langle\alpha^h(t_n^-),\beta^h(t_n^+)\rangle + \langle\rho c \vartheta^h(t_n^-), \sigma^{h}(t_n^+)\rangle+\langle\rho\Theta_0[\bm{\varepsilon}(\bm{u}(t_{n}^-))-\bm{\varepsilon}(\bm{u}(t_{n+1}^-))], \sigma^h(t_{n}^{-}) \rangle\\
&+(\bar{h},\sigma^h)_{F_n}+(\rho r, \sigma^h)_{Q_n}.
\end{aligned}
\end{equation}
The relation between the DG formulation \eqref{DG-T-Form} and the point-wise form \eqref{T-Strong} is apparent from the Euler-Lagrange form
\begin{equation}
\begin{aligned}
0&=A_{_n}^{^T}(\bm{\Pi}^h,\bm{\Lambda}^{h})-b_{_n}^{^T}(\bm{\Lambda}^{h})\\
&\quad + (\dot{\alpha}-\Theta^h, \beta^h)_{Q_n}\\
&\quad+(\rho \Theta_0\dot{\eta}+\mathrm{div}[\mathbf{k}_2\nabla\alpha^h+\mathbf{k}_3\nabla \theta^h]+\rho r,\sigma^h)_{Q_n}\quad\,\hspace*{5.7em} \text{(Equation of Motion)}\\
&\quad + \langle\jump{\alpha^h}_{n},\beta^h(t_n^+)\rangle\;\,\hspace*{18em}\text{($\alpha$-continuity)}\\
&\quad + \langle\jump{\rho c \vartheta^h}_n,\sigma^h(t_n^+)\rangle,\hspace*{17em}\;\text{($\vartheta$-continuity)}
\end{aligned}
\end{equation}
which reveals that a sufficiently smooth solution of the strong problem \eqref{T-Strong} also satisfies \eqref{DG-T-Form}, and vice versa, while the jump terms are vanished at the smooth solution. This also proves the consistency of the T-DG scheme of the thermal problem.

The unconditionally stability of the scheme \eqref{DG-T-Form} can also be shown along the same line of argument used for the mechanical case.
\subsubsection*{Remark:}
\begin{itemize}
\item[] Again, the use of the $L^2$-inner product to enforce the thermal problem allows one to omit the boundary restriction when we define the thermal displacement trial and weight function space. i.e.
\begin{equation}
\mathcal{T}^{^\alpha}_{_h}\;{:=}\;\mathcal{S}_{_{h}}(Q_n;j,l)\;{=:}\;\mathcal{W}^{^\alpha}_{_h}.
\end{equation}
This is a very important observation in terms of practical implementation.
\end{itemize}
As we have seen from Section~\ref{sec:alg-OS} that consistent and stable sub-algorithms render a consistent and stable global algorithm in the sense of time-stepping algorithms based on operator-splitting. Both the mechanical and the thermal phase algorithms are shown to be consistent and unconditionally stable. Therefore, the algorithm for the global problem based on Lie-Trotter-Kato product formula is consistent and unconditionally stable. Moreover, the convergence of the global scheme follows from the well known result stated below.

\begin{thm}[Lax Equivalence Theorem]
For \textbf{consistent} numerical approximations, \textbf{stability} is equivalent to \textbf{convergence}.
\end{thm}

\section{Numerical results}\label{sec:NumResults}
In this section, we present a range of results for type II and III problems of non-classical thermoelasticity. We start by comparing convergence of the proposed splitting scheme against a monolithic approach in which all the governing equations are discretised simultaneously using the time-DG finite element method. For this, a 1-D problem of non-dimensional form is considered. The result shows excellent agreement between the monolithic and the splitting scheme. Then we go on to present various results in 1-D and 2-D. The examples in this case are designed to illustrate two key features of the time-DG scheme: (i) its performance in solving problems that involves the propagation of sharp gradients without creating spurious oscillations; and (ii) its capability in capturing the unique aspects of non-classical theory, for example, propagation of thermal wave and its complex response due to the coupling of elasticity problem.     

The family of problems considered in this section are organized as follows. To analyse the rate of convergence and capability of the proposed scheme, non-dimensional form of a 1-D non-classical thermoelastic problem is presented in Section~\ref{sec:nonDim}. The performance of the splitting algorithm is examined in Section~\ref{sec:initProp} for an initial temperature pulse propagation in a two dimensional square plate under plane strain condition. Finally, in Section~\ref{sec:expCylin}, a quasi-static expansion of a thick walled, infinitely long cylinder in plane strain condition is presented, which is modelled as type I and type III theory of thermoelasticity and the remarkable difference of thermal responses between the two models are also analysed. 

\subsection{Non-dimensional 1-D GNT}\label{sec:nonDim} 
The non-dimensional form of 1-D GNT problem given in \eqref{coupled-hyperbolic1} is 
\begin{equation}\label{non-dimensional}
\begin{aligned}
\partial_{_\tau}\bar{u} &= \bar{v},\\
\partial_{_\tau}\bar{v} &= \partial_{_\xi}[\varepsilon_{_1}\partial_{_\xi}\bar{u}-\bar{\vartheta}]+\bar{b},\\
\partial_{_\tau}\bar{\alpha} &= \bar{\vartheta},\\
\partial_{_\tau}\bar{\vartheta} &= \partial_{_\xi}[\partial_{_\xi}\alpha+k\partial_{_\xi}\bar{\vartheta}]-\varepsilon_{_{2}}\partial_{_\xi}\bar{v} + \bar{s},\\
\end{aligned}
\end{equation} 
with the dimensionless parameters
\begin{equation}
\varepsilon_{_1} = (\dfrac{v_{_f}}{v_{_s}})^{2},\qquad \varepsilon_{_2}=\dfrac{\Theta_0 m^2E}{\rho c},\quad\text{ and }\quad k=\dfrac{k_3}{\sqrt{\rho c}},
\end{equation}
where $\varepsilon_{_1}$ denotes the square of ratio of uncoupled velocities of the mechanical wave (or first sound) and thermal wave (or second sound), $\varepsilon_{_2}$ denotes the strength of the thermomechanical coupling, $k$ represents the non-dimensional classical heat conductivity. The speed of first sound $v_f$ is actually the speed of sound in the medium, that is  
\begin{equation}
v_f = \sqrt{\dfrac{E}{\rho}},
\end{equation}
where $E$ denotes the Young's modulus of the medium, while that of the second sound $v_s$ is a characteristic feature of the theory of non-classical heat conduction by Green and Naghdi  that represents the speed in which a thermal disturbance travels through the medium: 
\begin{equation}
v_s = \sqrt{\dfrac{k_2}{\rho c}}.
\end{equation}
The non-dimensionless variables are given by
\begin{equation}
\xi = x_{_c}^{-1}x,\quad \tau = t_{_c}^{-1}t,\quad \bar{u} =u_{_c}^{-1}u,\quad \bar{v} = \dfrac{t_{_c}}{u_{_c}}v,\quad\bar{\alpha} = \alpha_{_c}^{-1}\alpha,\quad \bar{\vartheta} = \dfrac{t_{_c}}{\alpha_{_c}}\vartheta,
\end{equation}
where $x_{_c}$, $t_{_c}$, $u_{_c}$, $\alpha_{_c}$ are characteristic quantities having the same dimension as $x$, $t$, $u$, $\alpha$ respectively that can be chosen according to the relations
\begin{equation}
\dfrac{x_{_c}}{t_c} = v_{_s}, \quad \dfrac{u_{_{c}}}{\alpha_{_c}}=\dfrac{mv_{_f}}{\rho}.
\end{equation}
From the above equations we can observe that there are infinitely many ways of choosing the characteristics constants without changing the form of the system \eqref{non-dimensional}.

The nondimensional energy counterpart of \eqref{energy}, also referred to as the $H^1$-norm, is given by
\begin{equation}\label{nonD-H1norm}
\mathcal{E}(\bar{\bm{\chi}}) = \int_{0}^{\bar{L}}\bigg[\varepsilon_{_1}[\partial_{_\xi}\bar{u}]^2+\bar{v}^2+\frac{1}{\varepsilon_{_2}}[\partial_{_\xi}\bar{\alpha}]^2+\frac{1}{\varepsilon_{_2}}\bar{\vartheta}^2\bigg]\mathrm{d}\xi,
\end{equation}
and the $L^2$-norm 
\begin{equation}\label{nonD-L2norm}
\Vert\bar{\bm{\chi}}\Vert^2 = \int_{0}^{\bar{L}}\big[\bar{u}^2+\bar{v}^2+\bar{\alpha}^2+\bar{\vartheta}^2\big]\mathrm{d}\xi,
\end{equation}
where $\bar{\bm{\chi}}=(\bar{u},\bar{v},\bar{\alpha},\bar{\vartheta})^{T}$ is the state vector at a given time. 
\subsubsection{Convergence}\label{sec:convRates}
For the purpose of the convergence analysis an exact solution to problem \eqref{non-dimensional} is obtained in such a way that source terms $\bar{b}$ and $\bar{s}$ are suitably prescribed such that a given state vector $\bar{\bm{\chi}}=(\bar{u}, \bar{\alpha},\bar{\alpha},\bar{\vartheta})^T$ be an exact solution \cite{Khalmonova2008}. To this end let the source terms be
\begin{equation}
\begin{aligned}
\bar{b} &= \frac{\pi^2}{4}\big[(\varepsilon_{_1}-1)\sin(\pi \xi)\sin(\pi \tau)+ \cos(\pi \xi)\cos(\pi \tau)\big],\\
\bar{s} &= \frac{\pi^2}{4}\big[k\pi\sin(\pi \xi)\cos(\pi t)+\varepsilon_{_2}\cos(\pi \xi)\cos(\pi \tau)\big],
\end{aligned}
\end{equation} 
so that the exact solutions are  
\begin{equation}
\begin{aligned}
\bar{u}(\xi, \tau) &= \bar{\alpha} = \frac{1}{4}\sin(\pi \xi)\sin(\pi \tau),\\
\bar{v}(\xi, \tau) &= \bar{\vartheta} = \frac{\pi}{4}\sin(\pi \xi)\cos(\pi \tau),
\end{aligned}
\end{equation}
defined on the space-time domain $(\xi, \tau)\in[0,\bar{L}]\times[0,\bar{T}]$. For convergence analysis the values of the non-dimensional parameters at taken at $\varepsilon_{_1} = 4$, $\varepsilon_{_2} = 0.2$, $k = 0$. Such set of values represents a strongly coupled problem of two purely hyperbolic systems (type II thermoelasticity). The space-time domain corresponds to $\bar{L}=1$ and $\bar{T}=0.25$. Bilinear finite element functions, $Q1$, are used in each space-time slab with each element having an aspect ratio of one (i.e. $h=\Delta\xi=\Delta \tau$). 

\begin{figure}
\centering
\begin{tabular}{c}
\includegraphics[angle=0, width=8.0cm,height=6.5cm]{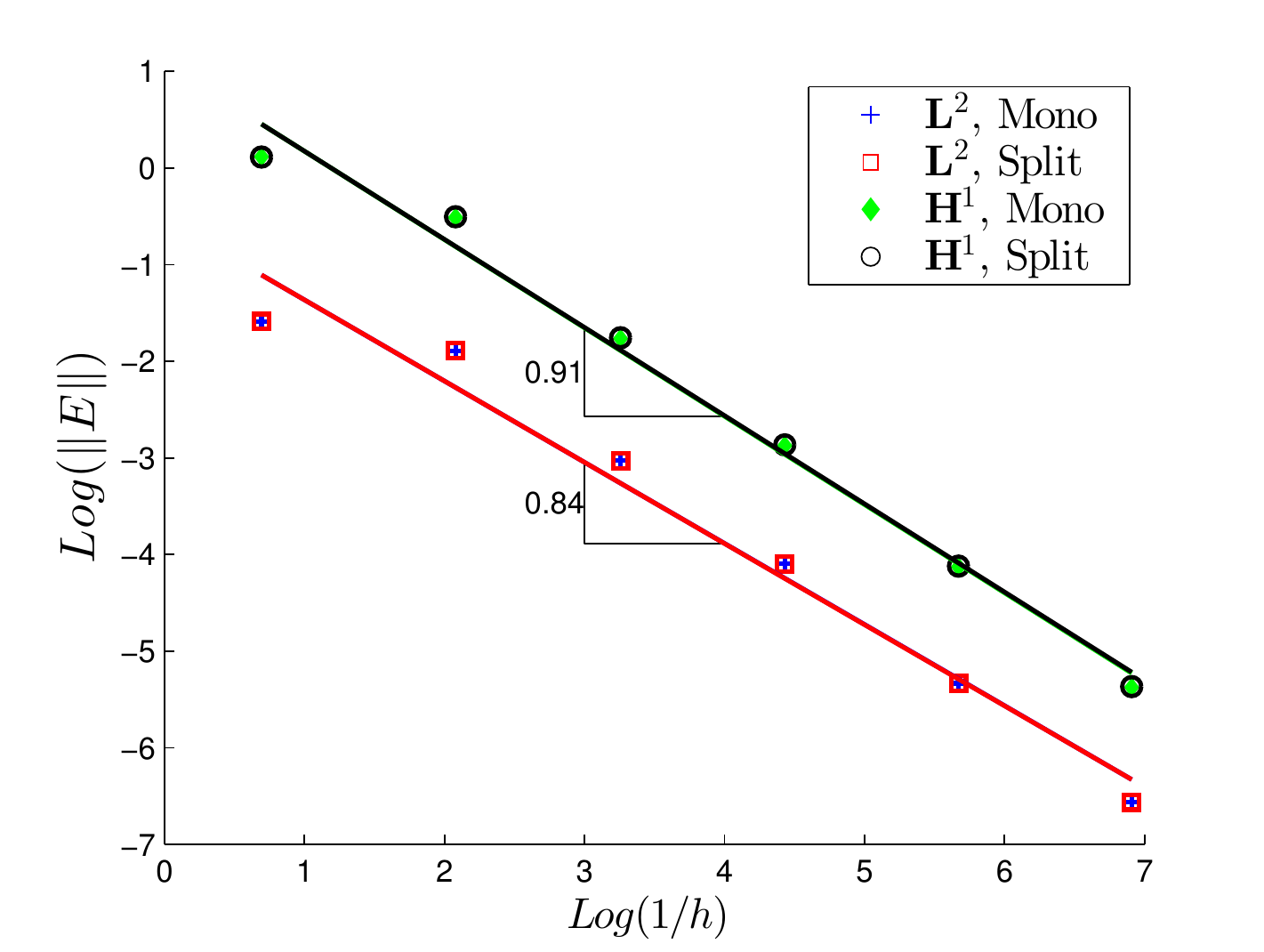}\\ 
 (a) Spatial convergence\\
\includegraphics[angle=0, width=8.0cm,height=6.5cm]{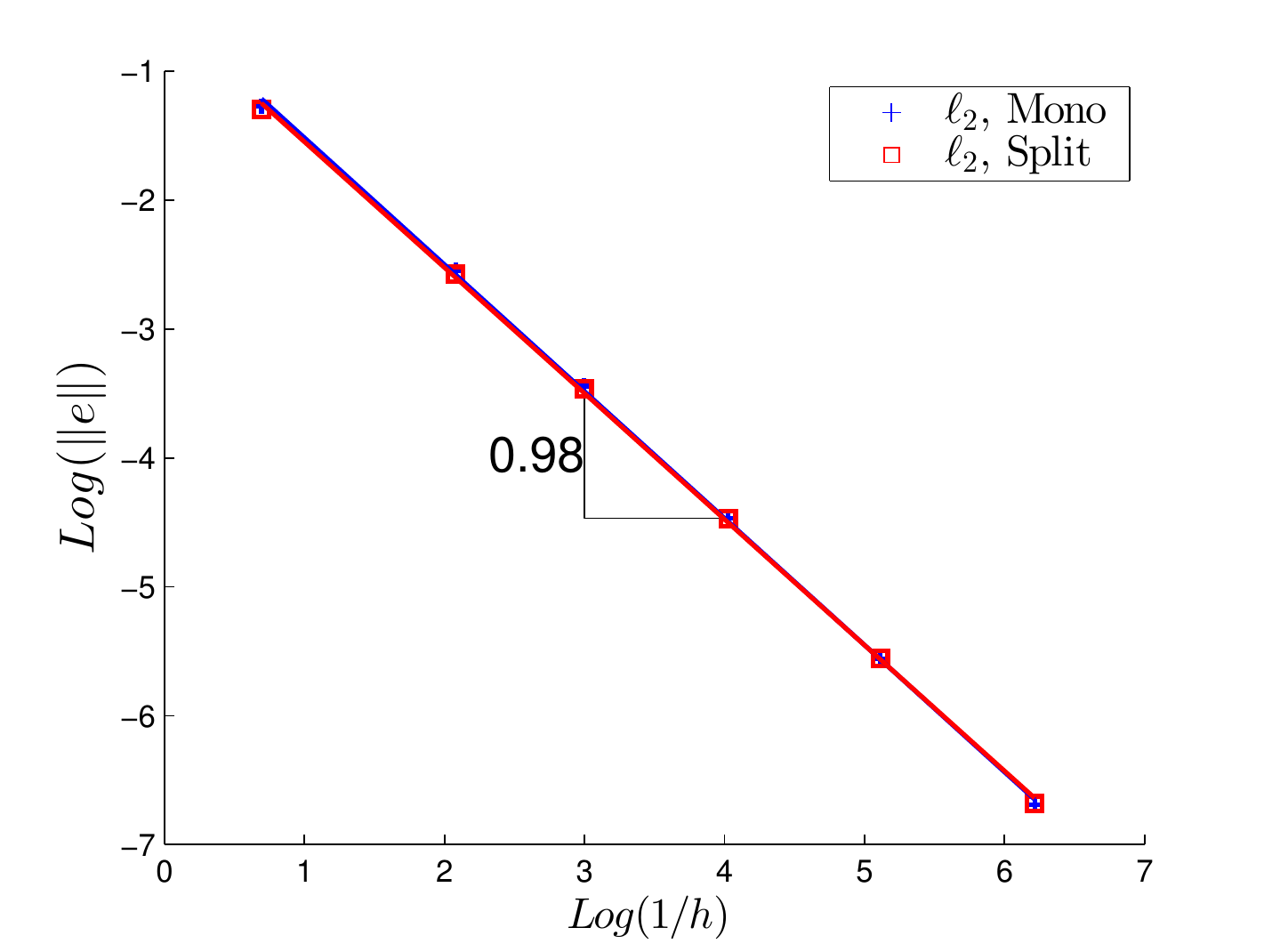}  \\ 
 (b) Temporal convergence
\end{tabular} 
\caption{Type II thermo-mechanical problem: Rate of convergence using monolithic and splitting approaches where the error norms are computed at $\tau = 0.25$ over the whole spatial domain.}
\label{fig:ConvRates}
\end{figure}

Fig.~\ref{fig:ConvRates}~(a) reports the \emph{spatial convergence} results of the monolithic and operator-splitting approaches with error norms of the approximate solutions at $\tau = \bar{T}$ are computed using the {$H^1$--} and $L^2$--norms as given in equations \eqref{nonD-H1norm} and \eqref{nonD-L2norm}, respectively . By construction the monolithic algorithm is only first-order accurate. Remarkably, it is shown that the error norms for both approaches are seen to overlap showing the error associated with the splitting is almost negligible. This demonstrably shows an increment in efficiency of the splitting algorithm, while maintaining the accuracy of the monolithic scheme.

While Fig.~\ref{fig:ConvRates}~(b) presents the result of \emph{temporal convergence} results of the two approaches with errors of the approximate solutions computed at the mid-point, $\xi=\bar{L}/2$, and $\tau=\bar{T}$ using the $\ell_{2}$ vector norm. Again errors corresponding to both of the approaches are almost identical. Similarly, this shows the temporal error associated with the operator-split strategy is minimal. 

\subsubsection{Laser pulse propagation}\label{sec:laserProp} 
Consider a one-dimensional bar occupying the interval $\xi\in[0, 1]$, heated by a pulsing laser applied at the left end having the form similar to the one considered in \cite{Miller2008} for non-Fourier heat conduction problem: 
\begin{equation}
\bar{s}(\xi, \tau) = \frac{1}{D\tau_{_p}}\exp\bigg[\bigg(\frac{\xi}{D}\bigg)^2-\bigg(\frac{\tau}{\tau_{_p}}\bigg)^2\bigg],  
\end{equation}
where $D$ is the depth of the pulse, and $\tau_{_p}$ is characteristic duration of the pulse. The bar is clamped at both ends at all times and with homogeneous initial conditions. We consider a situation in which a highly localized thermal pulse both in space and time described by the constants $\tau_{_p}=0.01$ and $D=0.02$ is applied at the left end of the bar. The parameters considered here are $\varepsilon_1=9$, which represents $3:1$ ratio of \emph{uncoupled} speeds of first sound to second sound,  and $\varepsilon_2=1$ accounting for a strongly coupled system. Bilinear elements are used in each space-time slab with mesh dimension $\Delta \xi = \Delta \tau = h$. The simulations are carried out over the period of $\bar{T}=1$ unit of non-dimensional time. The mesh parameter $h=0.001$ is chosen such that the width of the pulse is greater than the mesh size. In other words, the mesh is chosen so that the thin laser pulse can be described accurately by the the bilinear finite elements. 

\begin{figure}
\centering
\begin{tabular}{c}
\includegraphics[angle=0, width=12cm,height=6cm]{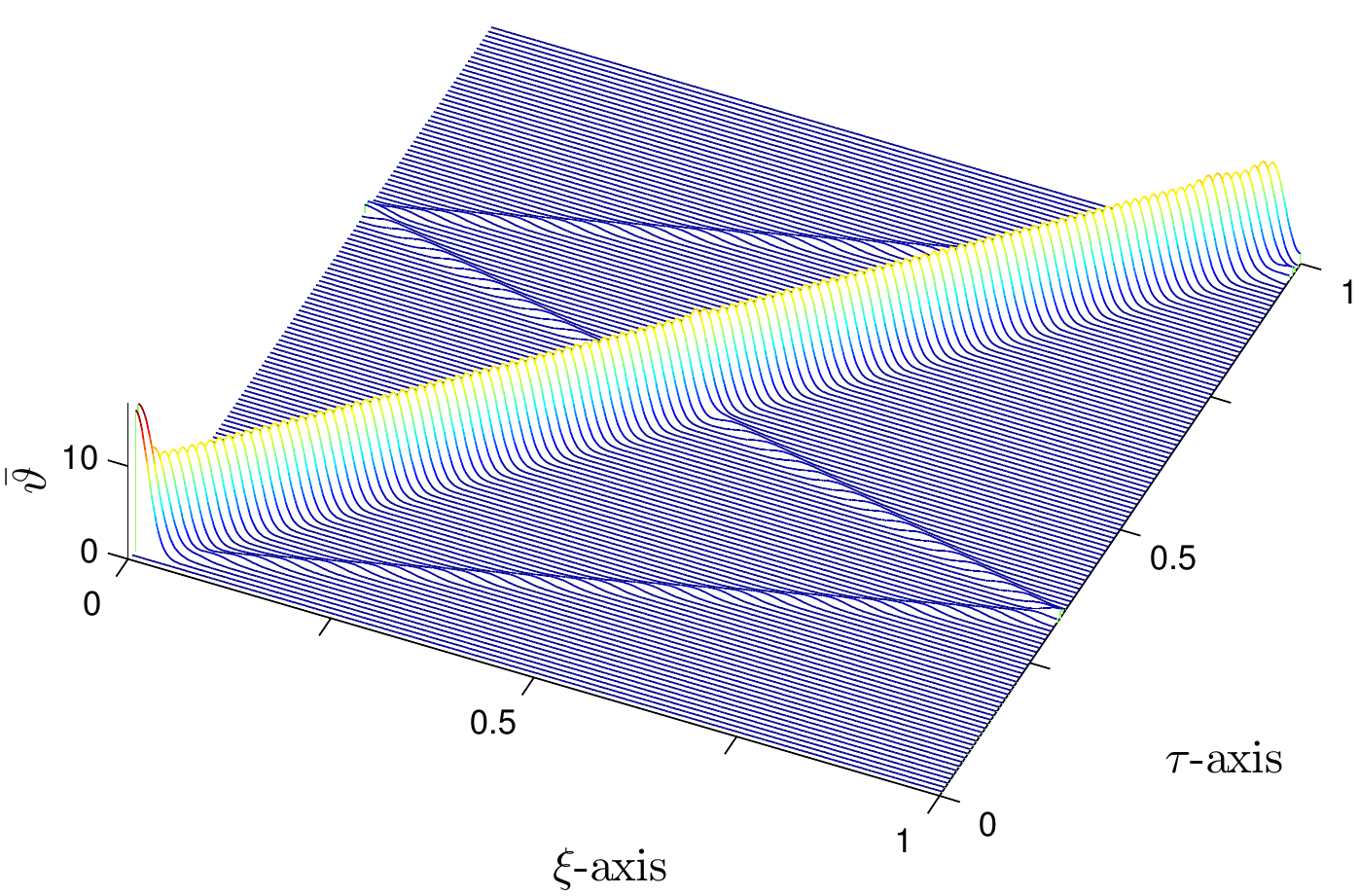}\\
(a) Monolithic\\ 
\includegraphics[angle=0, width=12cm,height=6cm]{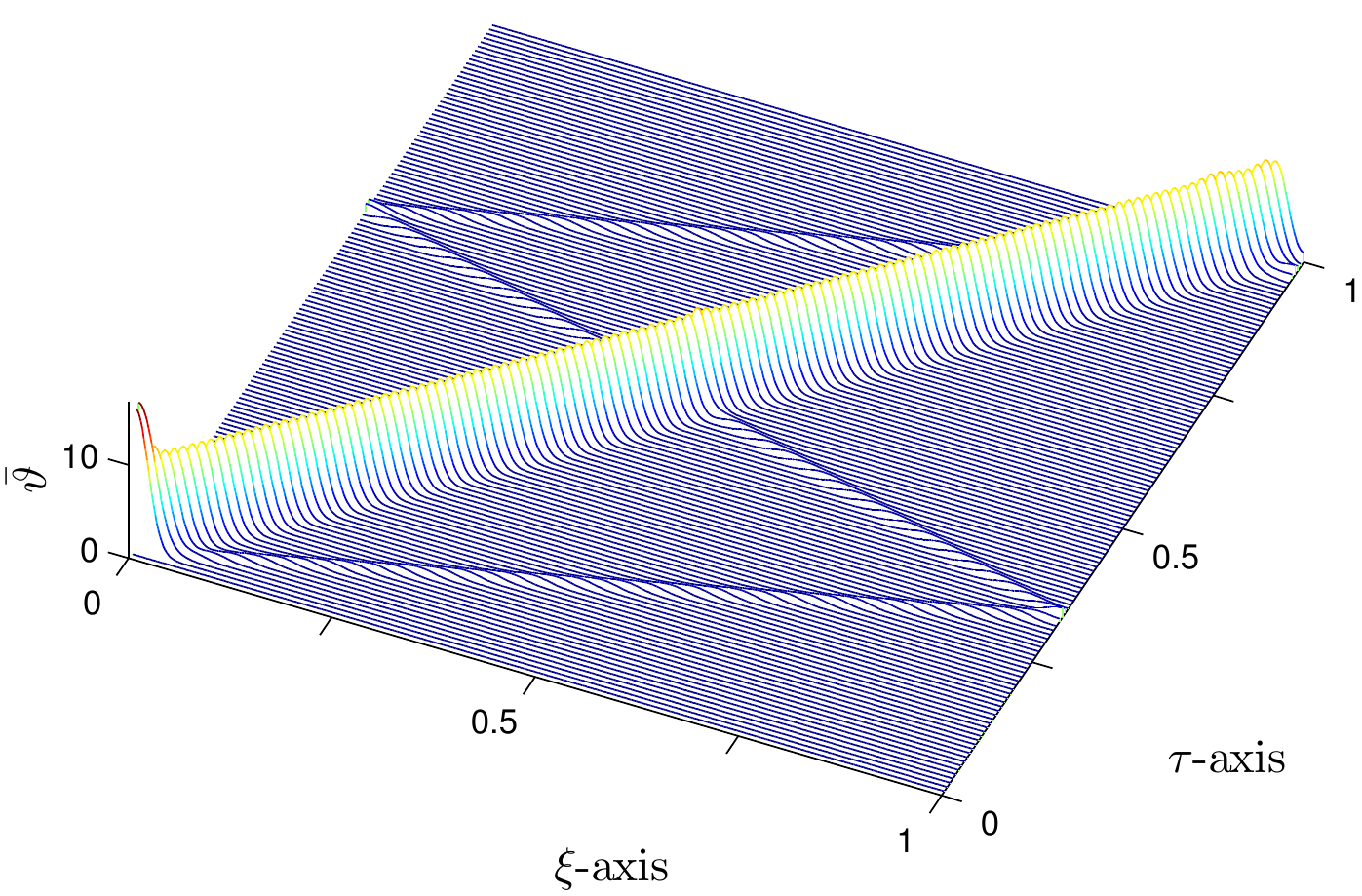}\\
(b) Splitting
\end{tabular} 
\caption{Propagation of laser pulse in type II thermoelasticity: temperature profile of the rod over the time period with $\varepsilon_{_{1}}=9$, $\varepsilon_{_{2}}=0.5$, $k=0$ and $\Delta \xi=\Delta \tau=0.001$.}
\label{fig:II_heatPulsing}
\end{figure}

Fig.~\ref{fig:II_heatPulsing} (a) and (b) show the propagation, in space and time, of the thermal disturbance caused by the pulsing laser heat source applied at the left end of the bar, computed using the monolithic and the splitting schemes, respectively. As can be seen from the figures, immediately after the pulse is applied, two thermal waves with different amplitude and speed start to emerge. The bigger and the slower wave is the one which is driven by the thermal equations, while the smaller and the faster one is induced by the mechanical equations through the coupling. The bigger thermal wave travels with a speed slightly less than that of second sound; whereas, the smaller thermal wave is travelling with a speed slightly greater than that of first sound. For this reason, it appears that the larger wave traverses the bar once, while the smaller traverses it more than three times. Note that the ratio of uncoupled speed of first to second sound is exactly $3:1$. 

Here there are two features which show the strength of the thermomechanical coupling: the first one is that the ratio of the speeds of the two thermal waves is noticeably different from what is expected in uncoupled case, and the other is that the coupling is strong enough to induce considerably large stress wave which in turn induce the faster thermal wave. 

Moreover, this problem represents a strongly coupled problem of two second-order hyperbolic problems involving propagation of sharp gradients. Such a problem is typically very difficult to approximate using the standard semi-discrete approach (MoL) unless some kind of stabilization term (or an artificial viscosity) is added , which is basically equivalent to changing the system from non-dissipative to dissipative, or very fine mesh is used together with very small time-step, which is undesirable from a computational cost point of view. 

What is remarkable about the current scheme is that it resolves the propagation of high gradients accurately while the amplitude of the thermal waves appear to be constant showing a very small numerical dissipation is added enough to damp out any numerical oscillation that could happen. The two approximate solution profiles Fig.~\ref{fig:II_heatPulsing} (a) and (b) are nearly identical. The agreement demonstrates that the splitting scheme maintains the accuracy of the monolithic scheme while the efficiency is considerably improved by the splitting scheme since two smaller systems are solved at each space-time slab. The result obtained here can be qualitatively compared to the one obtained in \cite{Bargmann2008}.

\begin{figure}
\centering
\begin{tabular}{c}
\includegraphics[angle=0, width=8.5cm,height=6.8cm]{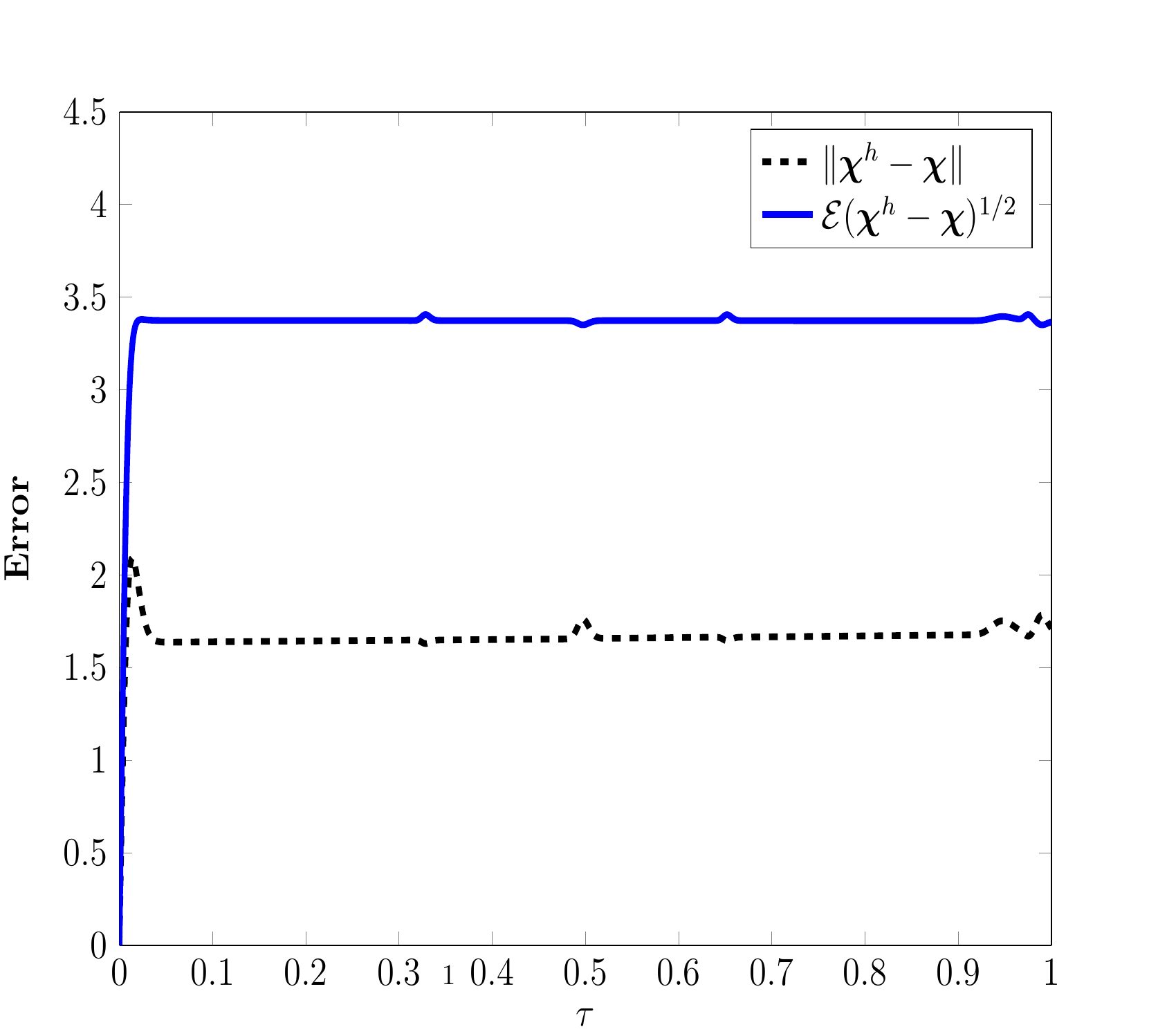}\\ 
(a) Monolithic\\
\includegraphics[angle=0, width=8.5cm,height=6.8cm]{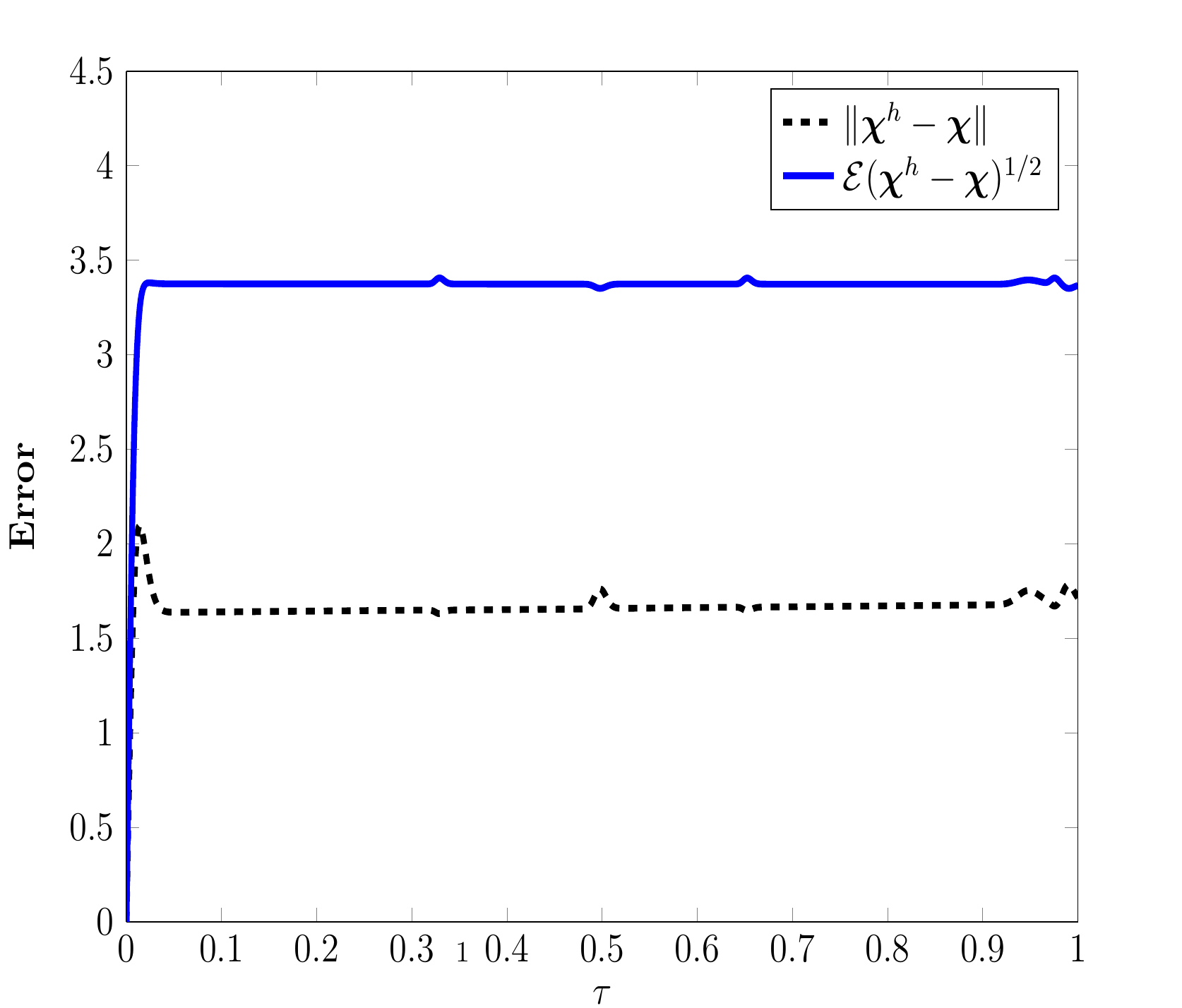} \\  
(b) Splitting
\end{tabular} 
\caption{Propagation of laser pulse in type II thermoelasticity: the $H^1$-Energies and $L^2$-norms corresponding to using monolithic and splitting approaches.}
\label{fig:Energies}
\end{figure}

As shown from Fig.~\ref{fig:Energies}, other than some small numerical instabilities when the waves interact either with the boundary or each other, the energy gained computed using the $H^1$-norm, remains essentially constant after the pulse is applied. This phenomenon is characteristic feature of type II thermoelasticity which is proved in Section~\ref{sec:well-posed}. While the $L^2$-norm shows more profound variation than the energy-norm right after the pulse is applied and when the two waves interact each other but it shows no change when the waves interact with the boundary. These observation suggests that the numerical instability that is arisen from the the interaction of waves with the boundaries may come from errors in the gradient of the approximate solution states. 

\begin{figure}
\centering
\begin{tabular}{c}
\includegraphics[angle=0, width=12cm,height=6cm]{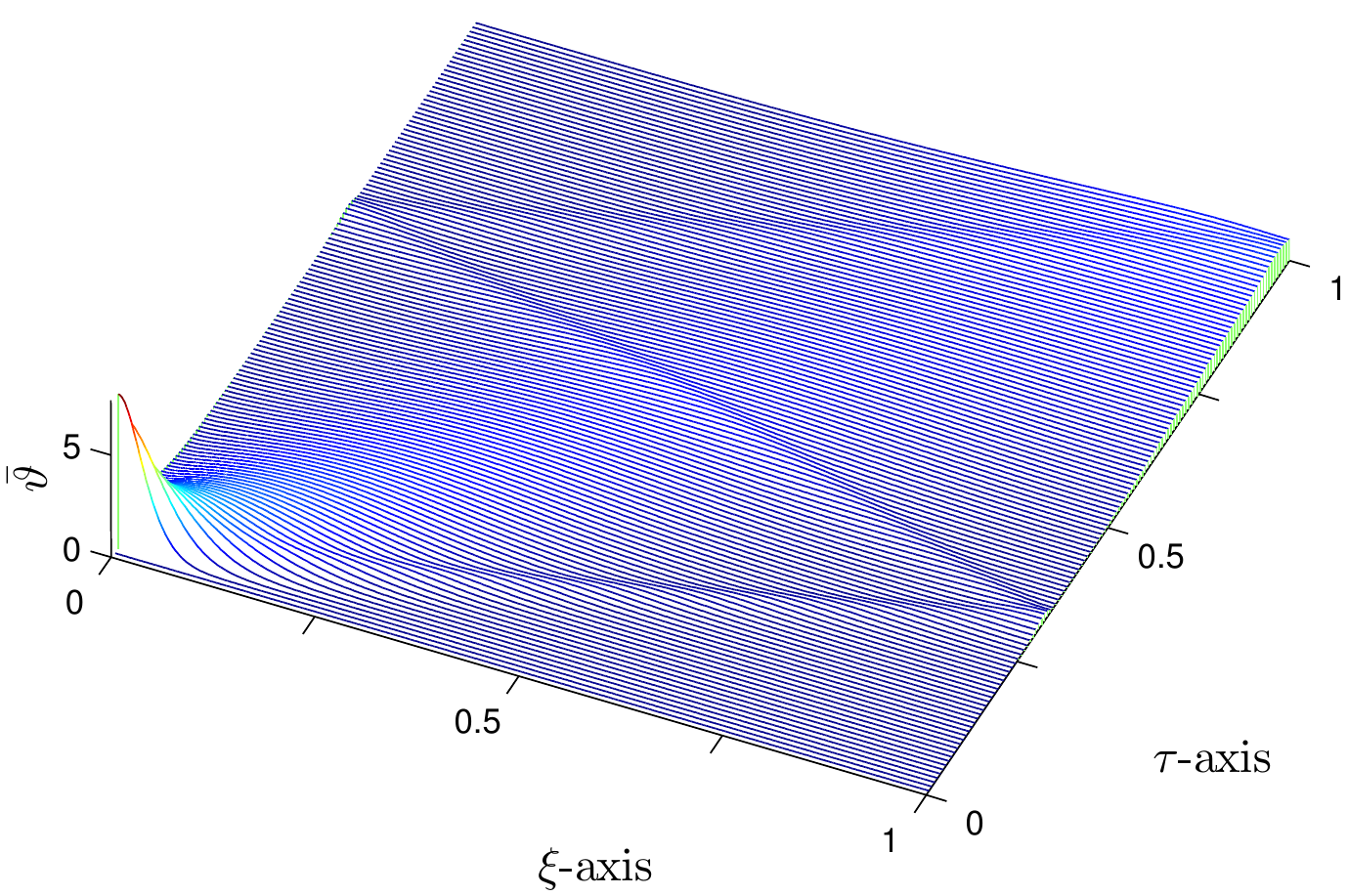}\\
(a) Monolithic\\
\includegraphics[angle=0, width=12cm,height=6cm]{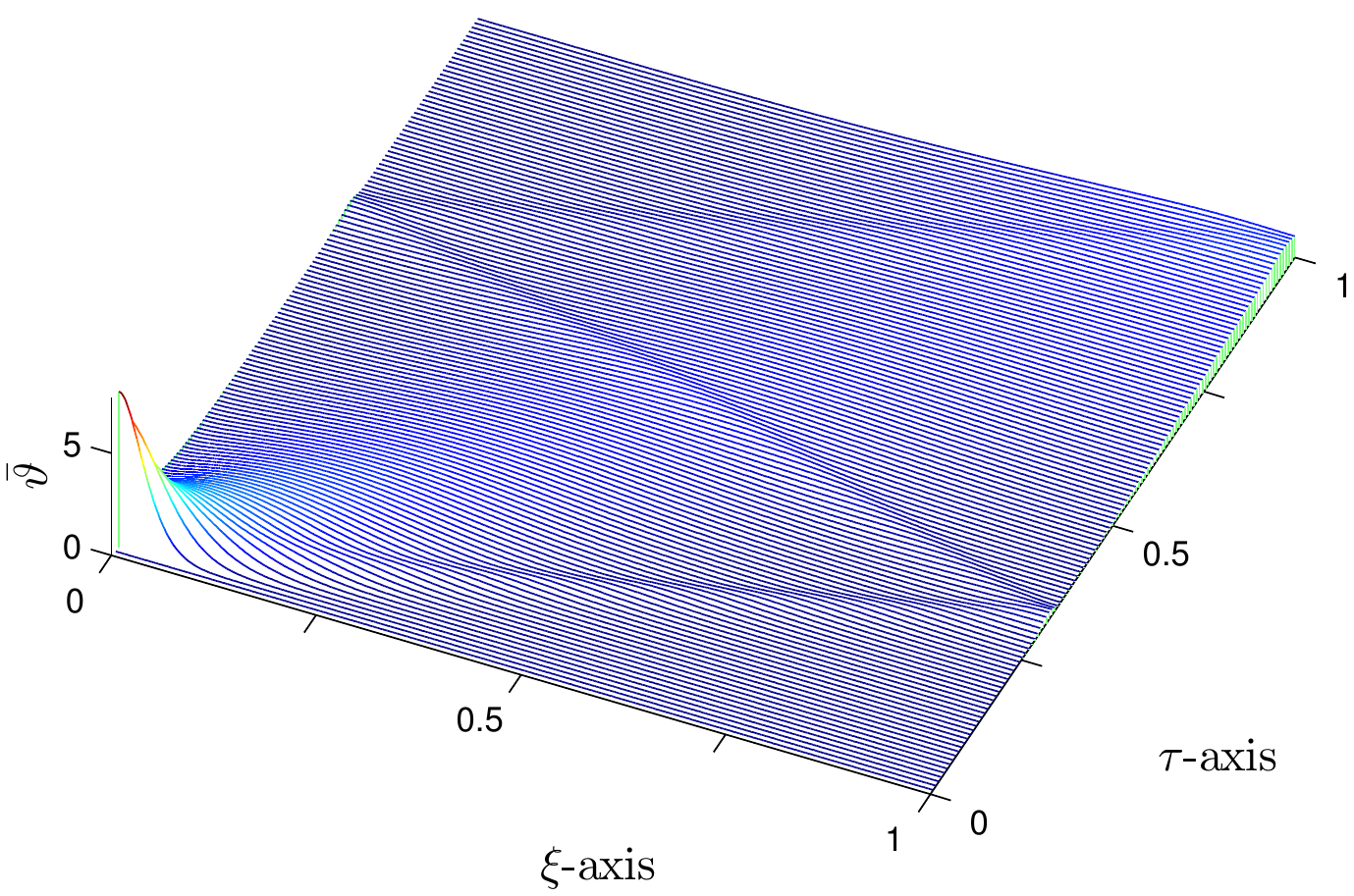}\\
(b) Splitting  
\end{tabular} 
\caption{Propagation of laser pulse in type III thermoelasticity: temperature profile of the rod over the time period with $k=0.1$, $\varepsilon_{_{1}}=9$, $\varepsilon_{_{2}}=0.5$, and $\Delta \xi=\Delta \tau=0.001$}
\label{fig:III_heatPulsing}
\end{figure}

Fig.~\ref{fig:III_heatPulsing}~(a) and (b) show the temperature profiles for the same problem above but with $k=0.1$, which correspond to Type III thermoelasticity, approximated using the monolithic and splitting schemes, respectively. This case is characterized by dissipation of energy while a wave scenario is still evident. The thermal wave driven by the temperature equations is damped out quickly, where as, the mechanically induced thermal wave remains localized for almost the entire duration  and is travelling with speed nearly equal to the speed of the first sound.  


\subsection{Two dimensional problem: Initial heat pulse propagation}\label{sec:initProp}
 In this problem, we consider a non-dimensional form of type III problem of initial thermal pulse propagation in a square plate occupying the region $\Omega=[-1, 1]\times[-1, 1]$ under plane strain assumption. A similar problem with the dimensions is solved in \cite{Bargmann2008}. The boundary of the specimen is mechanically clamped and fixed at the reference temperature $\Theta_0 = 1$ (i.e. the temperature of the ambient space). Initially, it was at rest but a temperature pulse is initialized at the center of the plate. i.e. the initial condition for the relative temperature $\vartheta$ be 
 \begin{equation}
 \vartheta(\bm{x}, 0) = A\exp\bigg[\dfrac{\bm{x}\cdot\bm{x}}{D}\bigg],
 \end{equation}
where $D$, as in the previous example in Section~\ref{sec:laserProp}, is a constant characterizing the width of the initial temperature pulse and $A$ is the amplitude. The material parameter used in the simulation are scaled according to the specifications summarized in Table~\ref{table:heatProp}.
\begin{figure}\label{tab:typeIII}
\begin{center}
\begin{tabular}{cc}
$t = 0$& $t = 0.2$\\
\includegraphics[angle=0, width=6.8cm,height=5.50cm]{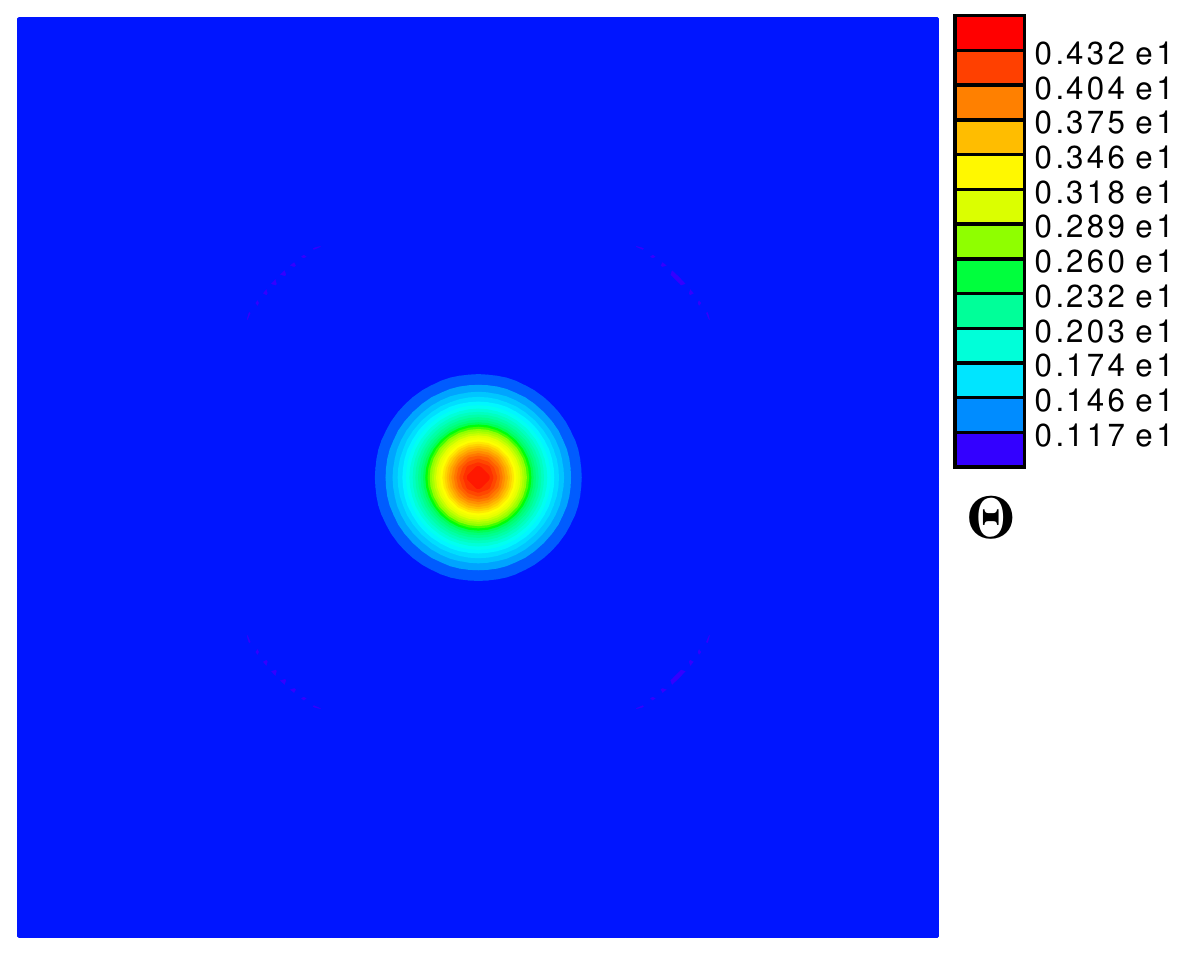}& 
\includegraphics[angle=0, width=6.8cm,height=5.50cm]{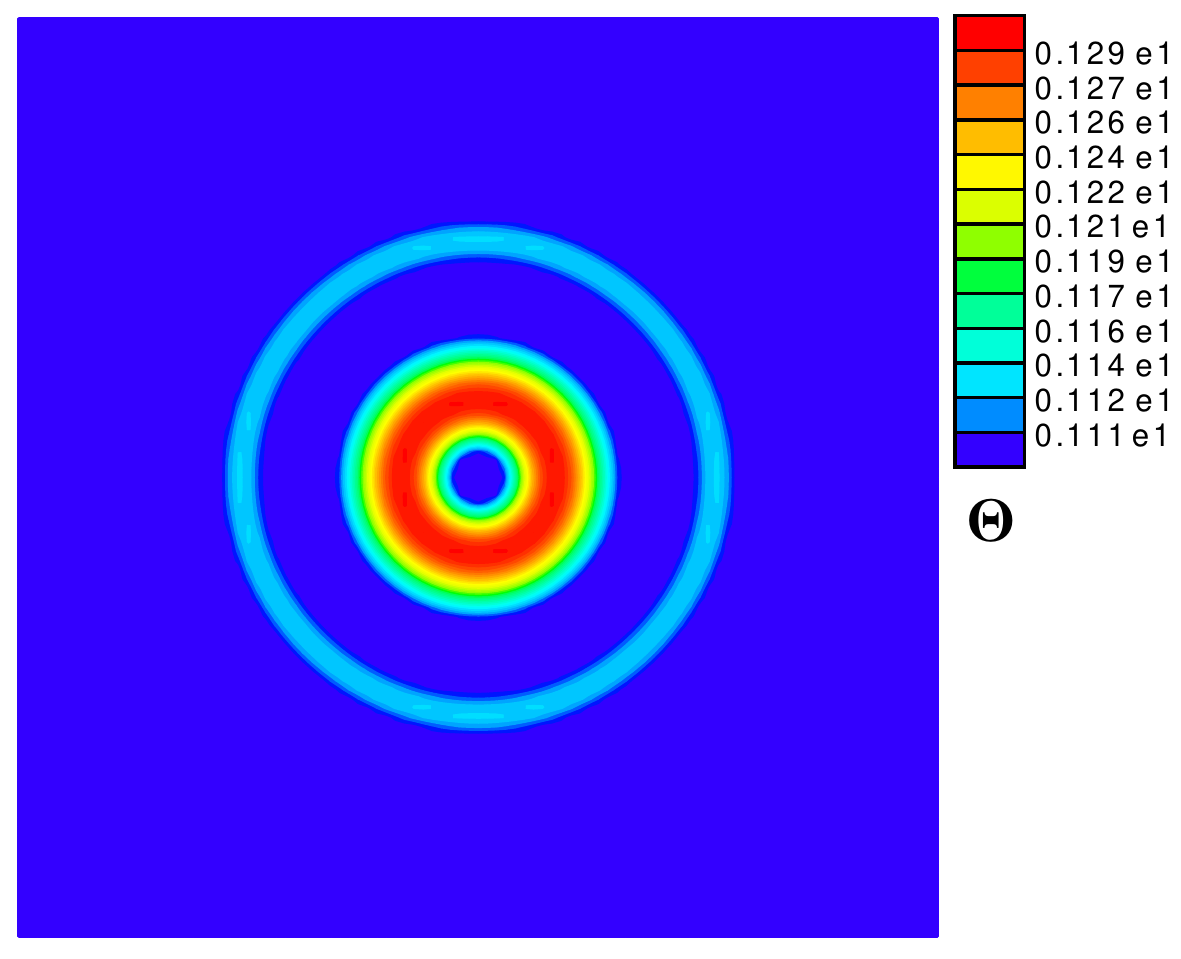}  \\ 
$t = 0.3$& $t = 0.4$\\
\includegraphics[angle=0, width=6.8cm,height=5.50cm]{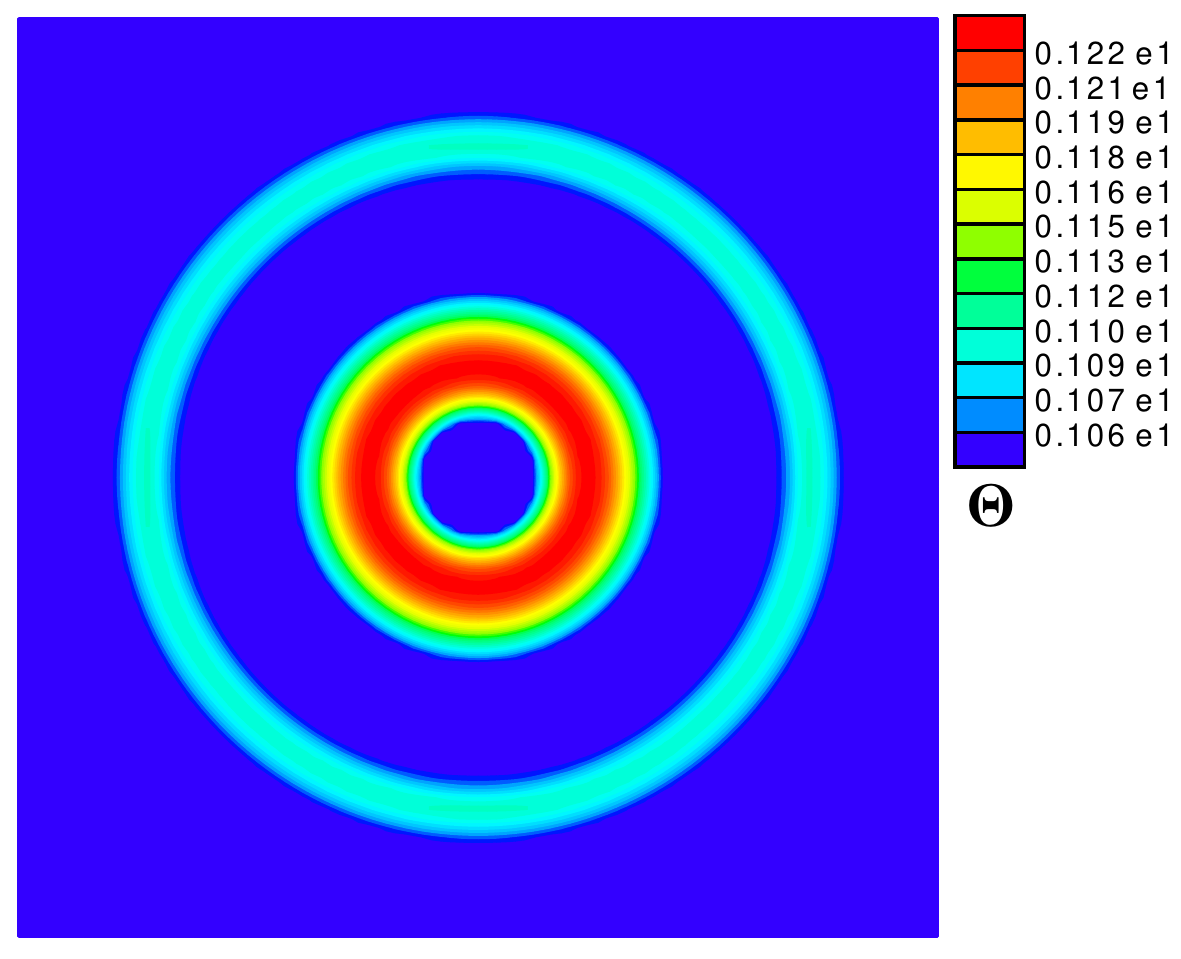}&
\includegraphics[angle=0, width=6.8cm,height=5.50cm]{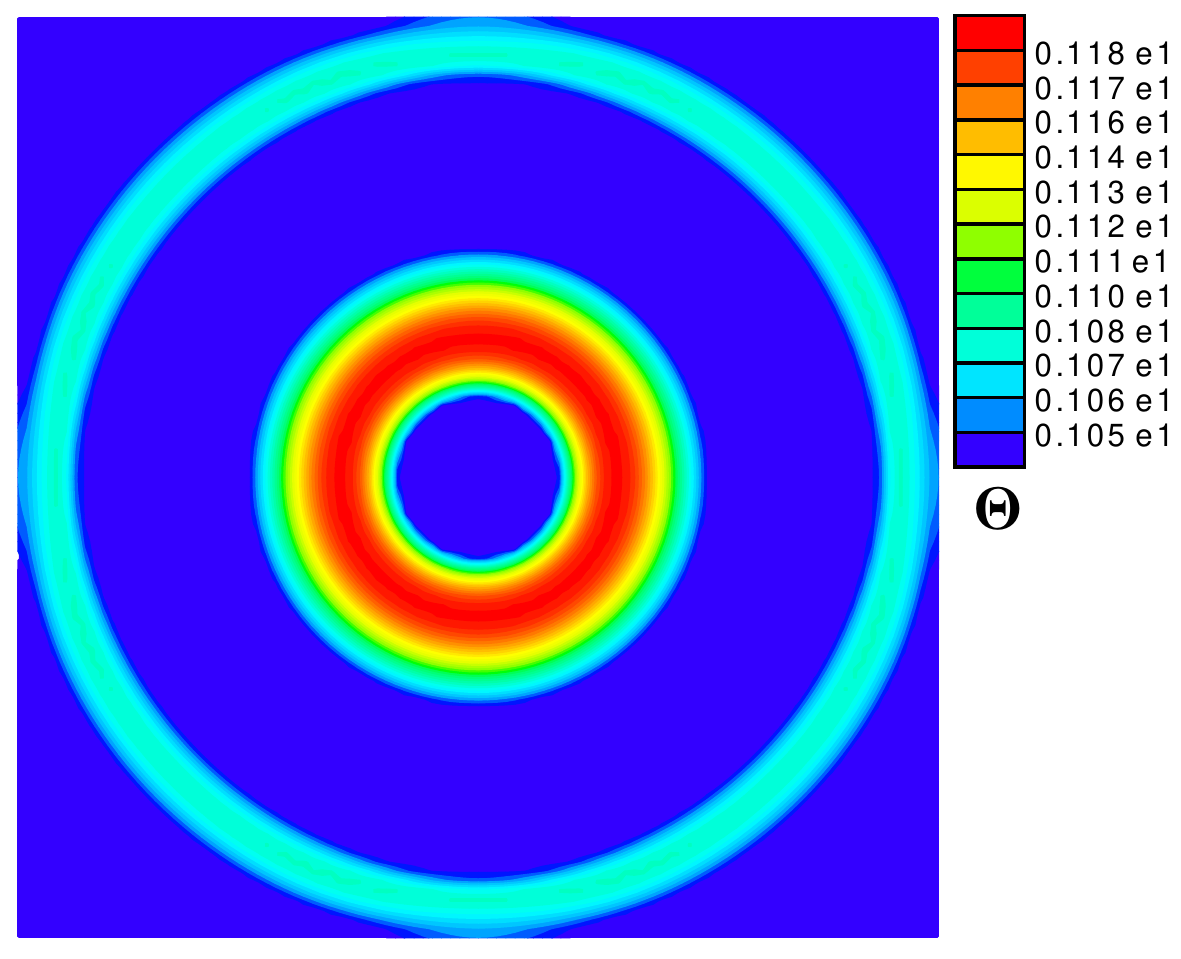}\\
\end{tabular} 
\caption{Temperature distribution in a square plate according to type III thermoelasticity where an initial pulse localized in space is initiated at the center.}
\label{fig:III_initPulsing}
\end{center}
\end{figure}
The time-DG finite element mesh consists of $8$ node isoparametric cubes with one element thickness $\Delta t = 0.01$ in time direction and $100\times100$ spatial elements per each slab are used to sufficiently describe the initial thermal pulse propagation.

Fig.~\ref{fig:III_initPulsing} shows snapshots of propagation of an initial temperature pulse with $D=100$ and $A=4$ at times $t=0$, $t=0.2$, $t=0.3$, and $t=0.4$. The initial pulse may be thought of as a thermal configuration just after an intense and highly localized laser heat source is applied at the center. The temperature profile gradually widens and a smaller but faster mechanically driven wave emerges, while the the second sound wave is driven by the temperature equations moves with a slower speed. In this case, the  classical conductivity parameter $\kappa_2$ gives additional stability but it is not so high to smear out the two wave phenomena.       

\begin{figure}
\begin{center}
\captionof{table}{Initial pulse propagation: material properties}  \label{table:heatProp}
\begin{tabular}{l c c}
\hline
Speed of first sound &  \hspace*{2em}$\sqrt{E/\rho}$ & \hspace*{2em}$1.96$\\
Speed of second sound &  \hspace*{2em}$\sqrt{\kappa_2/\rho c}$ &  \hspace*{2em}$0.65$\\
Conductivity ratio &  \hspace*{2em}$\kappa_2/\kappa_3$ &  \hspace*{2em}$100$\\
\hline
\end{tabular}
\end{center} 
\end{figure}

\subsection{Quasi-static case: Expansion of a thick walled cylinder}\label{sec:expCylin}
 This problem deals with the quasi-static thermo-mechanical interaction in a thick walled cylinder as it expands as a result of an inner wall Dirichlet-type boundary condition, in the linear and plane strain case. The material considered is isotropic both thermally and mechanically. The thermal variation is purely the result of mechanical changes (the expansion of the cylinder) unlike in the previous examples (Sections~\ref{sec:laserProp} and \ref{sec:initProp}) in which thermal variations cause mechanical effects. 
 
The cylinder has cross section occupying the region $\Omega = \{(x,y):r_0^2\leq x^2+y^2\leq R^2\}$ with inner and outer radii $r_0=10$ mm and $R = 20$ mm, respectively. A zero heat flux boundary condition is maintained on the inner wall, while the outer wall is kept at the reference temperature $\Theta_0$. The inner wall is dynamically prescribed a radial displacement of $1$ mm per second, while the outer wall is mechanically free.
 \begin{figure}
 \begin{center}
 \includegraphics[angle=0, width=5.0cm,height=5.0cm]{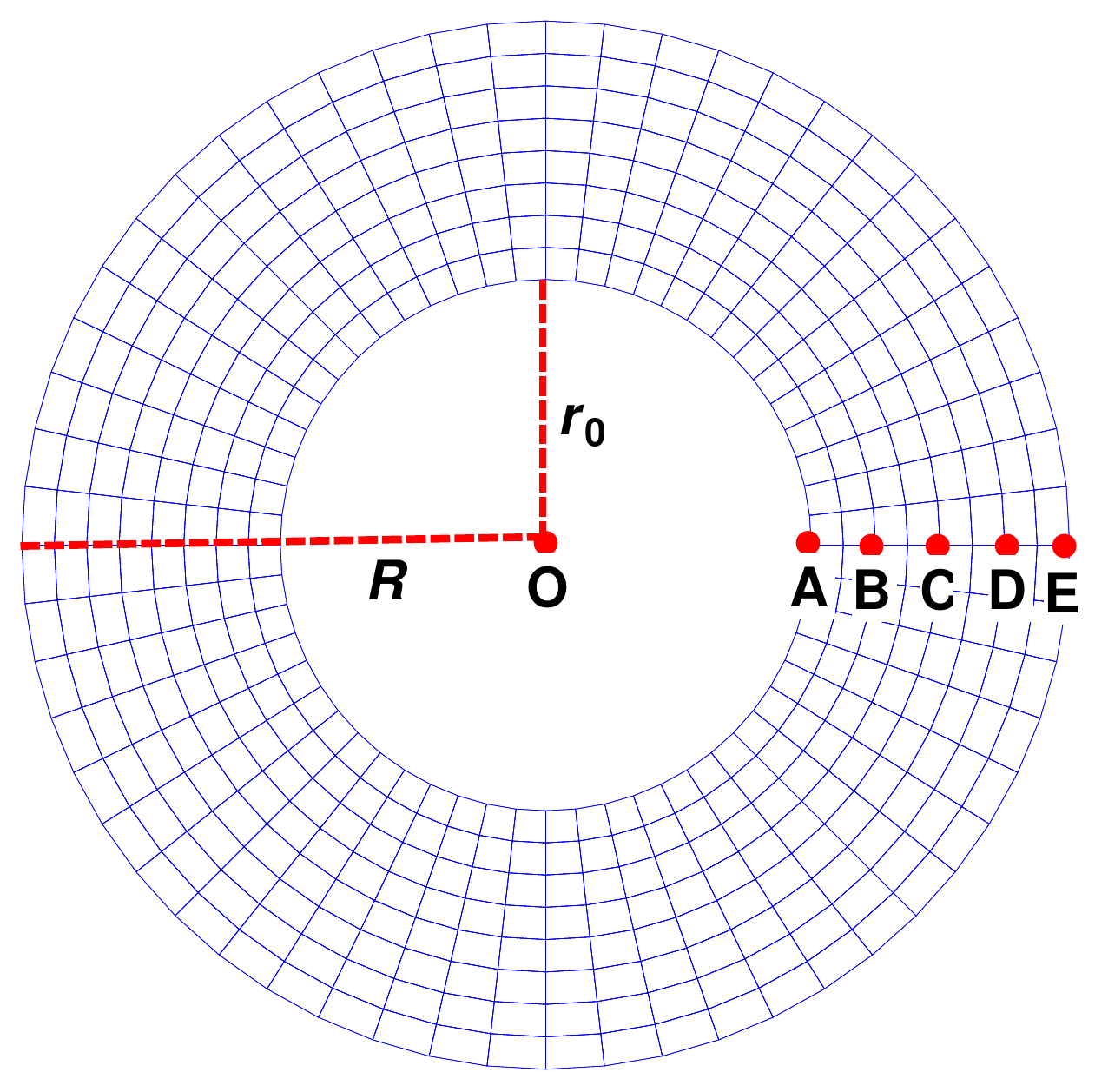}
 \end{center}
 \caption{Finite element mesh for the problem of expansion of a thick-walled cylinder.}
 \label{fig:cylinderMesh}
\end{figure}

  The problem is analysed for $20$ seconds until the inner wall reaches a radius of $r=3r_0$. The time-DG finite element mesh consists of trilinear shape functions of $56$ elements around the circumference of the cylinder by $8$ elements radially with one element thickness in the temporal direction with step length $\Delta t=0.1$ s for each space-time slab. In this quasi-static case, since only the thermal equations contains temporal derivatives then the thermal fields are allowed to be discontinuous while the displacement field is continuous across the interfaces of each space-time slab. This implies that the numerical dissipation comes from the weak enforcement of the continuity of the thermal fields only.

We consider two cases: the first is classical or type I thermoelasticity with  $k_3=45\text{ N/sK}$ and the 
other is type III thermoelasticity with $k_2 = 90\text{ N/K}$ and $k_3  = 30\text{ N/sK}$. 

\begin{figure}
 \centering
\begin{tabular}{c}

\includegraphics[angle=0, width=7cm,height=4cm,trim={0.0cm 0.0 0.0cm 0.0cm},clip]{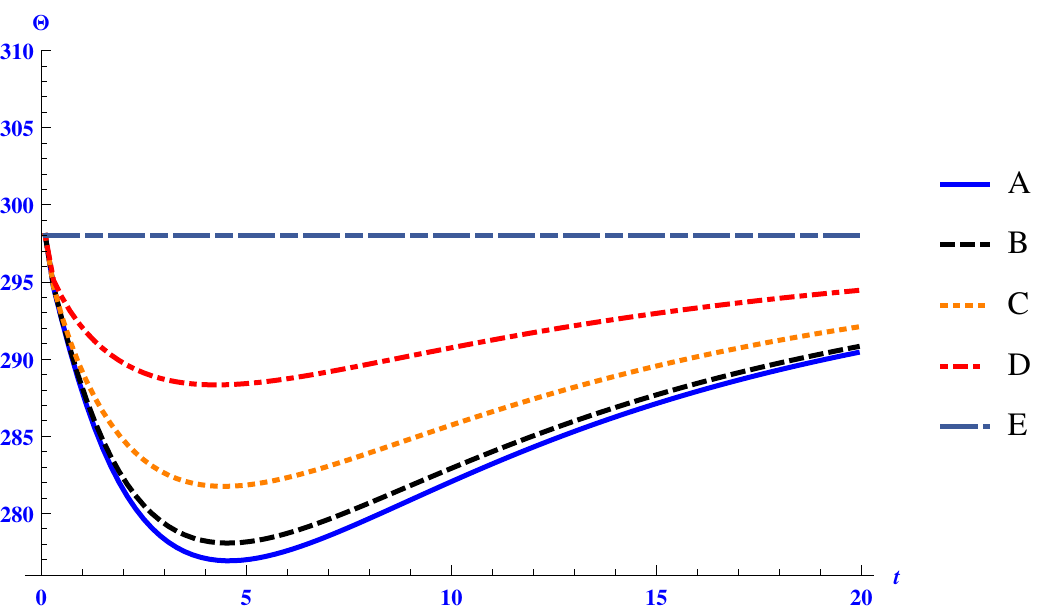}\\
(a) \textbf{Type I}\\
\includegraphics[angle=0, width=7cm,height=4cm,trim={0.0cm 0.0 0.0cm 0cm},clip]{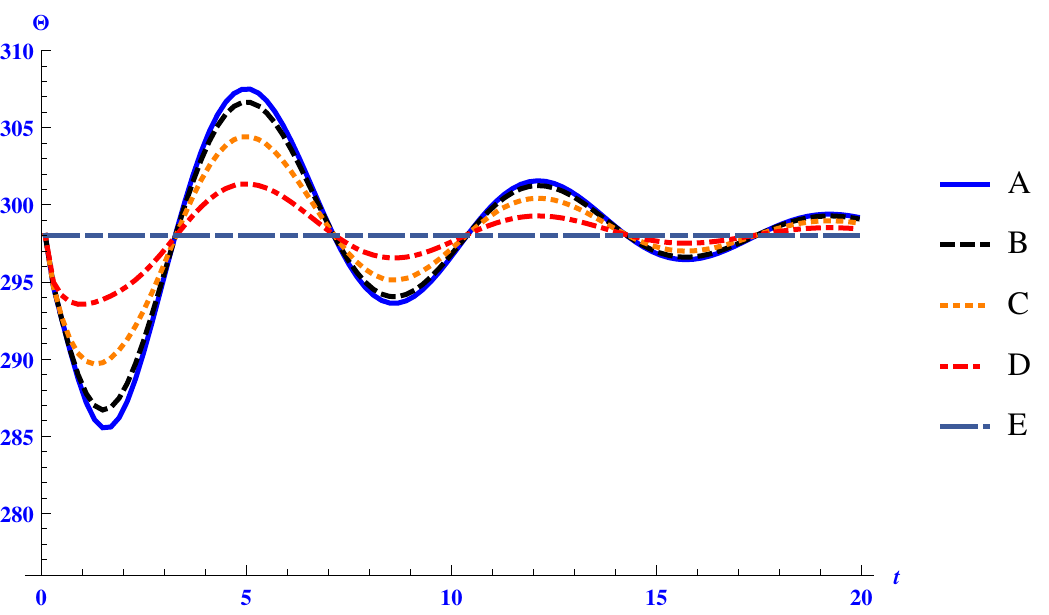}\\
(b) \textbf{Type III}
\end{tabular} 
\caption{Temperature profile of five points in the cylinder which are $0$ mm, $2.5$ mm, $5$ mm, $7.5$ mm, and $10$ mm away from the inner wall and shown with the labels A, B, C, D, and E.}
\label{fig:TemProfile}
\end{figure}

Fig.~\ref{fig:TemProfile} shows temperature variations over time for each case sampled at the equally spaced points along 
the radial direction labeled A-E as shown in the Fig.~\ref{fig:cylinderMesh}. As expected, in both cases, the 
temperature of the entire cylinder is converging to the reference temperature as time increases. The sinusoidal thermal response of type III is due to a temperature wave moving back and forth indicating second sound phenomenon.     

\section{Conclusion}\label{sec:conclusion}
 An operator-splitting strategy coupled with a space-time discontinuous Galerkin finite element method for the solution of transient and fully coupled initial-boundary problem of generalized thermoelasticity was presented. Well-posdeness of the problem in the general setting (type III) is proven using the theory of semigroups. The defining operator is split additively so that the first sub-operator represents an isentropic (adiabatic) elasticity in which
the entropy density is held fixed, and the other is a non-standard heat conduction at fixed configuration. Both of the sub-problems are also shown to inherit the same contractivity property as the full problem.

Each sub-problem is then discretised separately using a time-discontinuous Galerkin finite element method where the unknown fields are allowed to be discontinuous along the interfaces of each space-time slab. Weak continuity of the unknown fields is enforced using an $L^2$-inner product which differs from the original time-discontinuous formulation using an energy-inner product \cite{Hughes1988,Hulbert1990} which was formulated for linear elastodynamics problem.  The unconditionally stability behaviour of each of the algorithms is proven without the need to add extra `artificial viscosity'. The algorithm for the global problem is finally obtained by way of Lie-Trotter-Kato product formula, leading to an unconditional stability.     

The results presented in this paper are demonstrated by a number of numerical examples in both one and two dimensional cases. The efficiency of the current numerical scheme were examined by comparing the rate of convergence of with the corresponding monolithic approach. The result shows that the splitting scheme not only it retains the accuracy of the monolithic scheme but also it improves the efficiency as two smaller problems are solved sequentially at each time-step. The capability of the splitting algorithm is tested using problems involving propagation of heat waves driven by a pulsing laser heat source and an initial temperature disturbance in one and two dimensions respectively. Furthermore, the capability of the non-standard thermoelasticity and the proposed numerical method  to model the phenomenon of second sound in some solids is demonstrated by considering the quasi-static expansion of an infinitely long thick walled cylinder in plane stress.

The DG formulation proposed in this work may be extended to the non-linear regime
 without the need to eliminate the displacement-velocity relation in the formulation. Hence, a full recovery of the numerical dissipation in the non-linear case is possible. 
This will be the subject of a forthcoming work in which issues such as non-linear stability and the existence of Lyapunov function are discussed.

{\bf Acknowledgments.} The work reported in this paper has been supported by the National Research Foundation of South Africa through the South African Research Chair in Computational Mechanics. This support is acknowledged with thanks. The authors also thank Professor S. Bargmann for discussions which led to various improvements in the work.

\end{document}